\documentclass[final,12pt,a4paper]{article}

\usepackage{showkeys}
\usepackage[english,francais]{babel} 
\usepackage{ucs}
\usepackage[utf8]{inputenc}

\RequirePackage{color} 

\newif\ifpdf 
 \ifx\pdfoutput \undefined 
 \pdffalse 
\else 
 \pdfoutput=1 
 \pdftrue 
\fi

\ifpdf
 \usepackage[pdftex,colorlinks=true,urlcolor=blue,pdfstartview=FitH]{hyperref}
 \usepackage[pdftex]{graphicx}
 \DeclareGraphicsExtensions{.pdf, .png, .jpg} 
\else 
 \usepackage{hyperref}
 \usepackage{graphicx}
 \DeclareGraphicsExtensions{.eps, .ps} 
\fi

\RequirePackage{amsmath,amssymb,amsfonts,bbm,latexsym,mathrsfs}

\newcommand{\anc}{{\rm Anc}}
\renewcommand{\oe}{{\rm oe}}

\newcommand{\zzeta}{{\boldsymbol{\zeta}}}
\newcommand{\mmu}{{\boldsymbol{\mu}}}
\newcommand{\nnu}{{\boldsymbol{\nu}}}
\newcommand{\pphi}{{\boldsymbol{\varphi}}}

\newcommand{\dd}{{\partial}}

\newcommand{\bw}{{\bf w}}
\newcommand{\bW}{{\bf W}}

\newcommand{\ff}{{\bf f}}
\newcommand{\bs}{{\bf s}}
\newcommand{\bt}{{\bf t}}
\newcommand{\bz}{{\bf z}}

\newcommand{\bx}{{\bf x}}
\newcommand{\bX}{{\bf X}}

\newcommand{\bZ}{{\bf Z}}

\newcommand{\bQ}{{\bf Q}}

\newcommand{\bb}{{\bf b}}

\newcommand{\R}{\mathbb{R}}
\newcommand{\N}{\mathbb{N}}
\newcommand{\F}{\mathbb{F}}

\newcommand{\Z}{\mathbb{Z}}
\renewcommand{\P}{\mathbb{P}}

\newcommand{\E}{\mathbb{E}}

\newcommand{\D}{\mathbb{D}}
\newcommand{\T}{\mathbb{T}}

\newcommand{\TT} { {\cal T }}

\newcommand{\FF} { {\cal F }}

\def\build#1_#2^#3{\mathrel{
\mathop{\kern 0pt#1}\limits_{#2}^{#3}}}

\newcommand{\HT}{{\rm ht}}

\def\cq{$\hfill \square$}

\def\d{{\rm d}}
\def\a{{\cal A}}

\def\u{{\cal U}}

\def\W{{\cal W}}

\def\eps{\varepsilon}

\newcommand{\ba}{{\bf a}}
\def\ov{\overline}
\def\wh{\widehat}
\def\wt{\widetilde}

\def\rem{\noindent{\bf Remark. }}

\def\proof{\noindent{\bf Proof. }}

\def\card{{\rm Card\,}}

\def\var{{\rm Var\,}}
\def\cov{{\rm Cov\,}}

\newcommand{\ind}{\mathbbm{1}}

\newtheorem{thm}{Theorem}
\newtheorem{lmm}{Lemma}
\newtheorem{prp}{Proposition}
\newtheorem{defn}{Definition}
\newtheorem{crl}{Corollary}

\usepackage{fancyhdr}

\usepackage[outerbars]{changebar}

\def\miermont{{\href{http://mahery.math.u-psud.fr/~miermont/}{Gr\'egory
      Miermont}}}
\def\orsay{{\href{http://www.math.u-psud.fr/}{Laboratoire de
      Math\'ematique}}} 
\def\cnrs{{\href{http://www.cnrs.fr/}{CNRS}}}

\title{Invariance principles for spatial
  multitype Galton-Watson trees} 

\author{
  \miermont \thanks{\cnrs\ \&\ \orsay,
    \'Equipe Probabilit\'es, Statistique et Mod\'elisation, B\^at.\
    425, Universit\'e Paris-Sud, 91405 Orsay, France. {\tt
      Gregory.Miermont@math.u-psud.fr}}}

\begin{document}

\selectlanguage{english}

\maketitle




\begin{abstract}
  We prove that critical multitype Galton-Watson trees converge after
  rescaling to the Brownian continuum random tree, under the
  hypothesis that the offspring distribution has finite covariance
  matrices. Our study relies on an ancestral decomposition for marked
  multitype trees. We then couple the genealogical structure with a
  spatial motion, whose step distribution may depend on the structure
  of the tree in a local way, and show that the resulting discrete
  spatial trees converge once suitably rescaled to the Brownian snake,
  under some suitable moment assumptions.
\end{abstract}

\bigskip

\noindent{\bf M.S.C. }60J80, 60F17

\bigskip

\noindent{\bf Key Words: }Multitype Galton-Watson tree, discrete
snake, invariance principle, Brownian tree, Brownian snake

\newpage

\section{Introduction and main results}\label{sec:intr-main-results}

\subsection{Motivation}

Multitype Galton-Watson (GW) processes arise as a natural generalization of
usual GW processes, in which individuals are differentiated by
{\em types} that determine their offspring distribution. They were
first studied in 1947 by Kolmogorov and his coauthors. We refer to
\cite[Chapters 2 \& 3]{harris63} for a very nice introduction to these
processes. It turns out that their analysis is considerably eased
under an {\em irreducibility} assumption, namely, that every type has
a positive probability to eventually `lead' to all others.  Under this
hypothesis, one can use variants of the Perron-Frobenius theorem,
which allow to quantify the asymptotic behavior of iterates of the
{\em mean operator} of \cite{harris63}, and obtain qualitative and
quantitative results on the GW process. Informally, the
large-scale aspects of irreducible multitype GW processes
are similar to that of monotype GW processes whose mean
offspring distribution is the Perron eigenvalue of the mean
operator. On the other hand, Janson \cite{jansonurn05} has shown that
if the irreducibility assumption fails to hold, many different
behaviors can occur.

The aim of the present paper is to investigate the ancestor trees and
forests associated with irreducible GW processes, when the
total number of types is finite. Under criticality hypotheses on the
mean matrix, and a finiteness hypothesis on the covariance matrices of
the offspring distributions, we show that the height process of these
forests converges to a reflected Brownian motion, and hence behaves
asymptotically in a similar way as monotype processes (Theorem
\ref{T1}). Similar results are proved for family trees conditioned on
the number of their individuals, under extra exponential moments
assumptions (Theorem \ref{sec:convergence-gw-trees}). Although we do
not focus on a continuum tree formalism here, this says roughly that
under these hypotheses, multitype GW trees conditioned to
have $n$ individuals converge once suitably renormalized to the
Brownian continuum random tree \cite{aldouscrt93}.

We also show that these ancestor trees and forests, when coupled with
a spatial motion whose step distribution may locally depend on the
tree structure, converge to the so-called Brownian snake
\cite{legall99}, under mild extra hypotheses on the spatial motion step
distribution (Theorem \ref{T2}).

Applications of these results to random planar maps are discussed in
the companion paper \cite{mierinv}. The approach of the present paper
follows quite closely that of \cite{jfmgm05}.

\subsection{Multitype Galton-Watson processes}\label{sec:mult-galt-wats}

Let $K\in\N=\{1,2,\ldots\}$ be a positive integer.  Write
$[K]=\{1,2,\ldots,K\}$, and identify $\Z^{[K]},\R^{[K]}$ with
$\Z^K,\R^K$.  Suppose given distributions $(\mu^{(i)}, i\in[K])$ on
the space $\Z_+^K$ of integer-valued non-negative sequences of length
$K$. We will often use the notation $\mmu$ as a shorthand for the
$\R^K$-valued measure $(\mu^{(i)},i\in [K])$. 

A $K$-multitype
GW process with offspring distributions $\mmu$ is a
$\Z^K_+$-valued Markov process 
$$\bZ_n=(Z_n(i),i\in [K])\, ,\qquad n\geq 0$$
such that the law of $\bZ_{n+1}$ given $\bZ_n=\bz=(z_i,i\in[K])$ is
the same as $\sum_{i\in[K]}\sum_{l=1}^{z_i}\bX^{(i)}(l)$, where the
vectors $\bX^{(i)}(l), i\in[K],l\geq 1$ are all independent, and
$\bX^{(i)}(l)$ has law $\mu^{(i)}$ for all $l\geq 1$.

Otherwise said, this process can be considered as a model for
population evolution where each individual is given a type in $[K]$,
and where each type-$i$ individual gives birth to a set of individuals
with law $\mu^{(i)}$, this independently over individuals (although
the different components of $\bX^{(i)}$ may well be dependent).

We say that the process (or the measure $\mmu$) is {\em
  non-degenerate} if there exists at least one $i\in[K]$ so that
$\mu^{(i)}(\{\bz:\sum_j z_j\neq 1\})>0$.  Failure of this last
assumption entails that all particles a.s.\ give birth to exactly one
particle, and the study of the process boils down to that of a Markov
chain with values in $[K]$. All the processes that we consider here
are assumed to be non-degenerate.

For $i,j\in[K]$, let
$$m_{ij}=\sum_{\bz\in\Z_+^{K}}z_j\mu^{(i)}(\{\bz\})
$$
be the mean number of type-$j$ offspring of an type-$i$ individual.
We let $M_\mmu=(m_{ij})_{i,j\in[K]}$ and call it the {\em mean matrix}
of $\mmu$.

\begin{defn}\label{irred}
  The mean matrix (or the offspring distribution $\mmu$) is called
  {\em irreducible}, if for every $ i,j\in[K]$, there is some $n\in\N$
  so that $m^{(n)}_{ij}>0$, where $m^{(n)}_{ij}$ is the $ij$-entry of
  the matrix $M_\mmu^n$.
\end{defn}

Notice that the $n$ in the definition may depend on the choice of
$i,j$, so that the definition is distinct of that of {\em
  aperiodicity}, namely that the above property holds jointly for
every $i,j$, for some common $n$.

With the irreducibility assumption, the Perron-Frobenius Theorem
recalled below (Proposition \ref{frobenius}) ensures that the
eigenvalue $\varrho$ of $M_\mmu$ with maximal modulus is real,
positive, and simple, and that a non-zero eigenvector of $M_\mmu$ with
eigenvalue $\varrho$ has only non-zero entries all having the same
sign.

\begin{prp}[\cite{AthNey}]\label{sec:mult-galt-wats-2}
  Suppose that the process $(\bZ_n,n\geq 0)$ is non-degenerate. Then
  it eventually becomes extinct a.s., whatever the starting value, if
  and only if $\varrho\leq 1$. The process (or the distribution
  $\mmu$) is called sub-critical if $\varrho<1$, critical if
  $\varrho=1$ and supercritical if $\varrho>1$.
\end{prp}

It will be useful to introduce the generating functions 
$\pphi=(\varphi^{(i)},i\in[K])$ defined by 
$$\varphi^{(i)}(\bs)=
\sum_{\bz\in\Z_+^{K}}\mu^{(i)}(\{\bz\})\bs^{\bz}$$
where $\bs=(s_i,i\in[K])\in[0,1]^K$ and $\bs^{\bz}=\prod_{i\in[K]}s_i^{z_i}$. 
With these notations, we have 
$m_{ij}=\partial \varphi^{(i)}/\partial s_j({\bf 1})$, where 
${\bf 1}$ is the vector of $[0,1]^K$ with all components equal to $1$. 

For $1\leq i,j,k\leq K$, define
$$Q^{(i)}_{jk}=
\frac{\partial^2 \varphi^{(i)}}{\partial s_j\partial s_k}({\bf 1}).$$
We say that $\mmu$ has {\em finite variance} if 
\begin{equation}\label{eq:6}
  Q^{(i)}_{jk}<\infty \, \quad \mbox{ for all }\quad i,j,k\in[K]\, .  
\end{equation}
Under this assumption, for each $i$, $(Q^{(i)}_{jk},1\leq j,k\leq K)$
is the Hessian matrix of the convex function $\varphi^{(i)}$ evaluated
at ${\bf 1}$, hence the matrix of a non-negative quadratic form on
$\R^K$, which we call $Q^{(i)}(\bs),\bs\in\R^K$.

Assuming $M_\mmu$ irreducible and $\mmu$ non-degenerate, critical and
with finite variance, we let $\ba,\bb$ be the left and right
eigenvectors of $M_\mmu$ with eigenvalue $1$, chosen so that $\ba\cdot
{\bf 1}=1$ and $\ba\cdot \bb=1$, where $x\cdot y$ is the scalar
product of the vectors $x,y\in\R^K$. Let
\begin{equation}\label{sigma}
\sigma=\sqrt{\sum_{i=1}^K a_{i}Q^{(i)}(\bb)}=\sqrt{\ba\cdot \bQ(\bb)}, 
\end{equation}
where $\bQ(\bs)$ is the $K$-dimensional vector $(Q^{(i)}(\bs),1\leq
i\leq K)$.  This should be interpreted as the `variance' of the
offspring distribution of the multitype process, as it plays a role
similar to the variance for monotype GW processes in the asymptotics
of the survival probability, see \cite{VD}.

\bigskip

\noindent{\bf Basic assumption (H). }
In the sequel, unless specified otherwise, we will exclusively be
concerned with irreducible, non-degenerate, critical offspring
ditributions with finite variance. Notice that criticality implies
finiteness of all coefficients of $M_\mmu$.

\subsection{Multitype trees and forests}\label{sec:mult-galt-wats-1}

We now add a genealogical structure to the branching processes, 
by endowing it into a tree-valued random variable. For $n\geq 0$
let $\u$ be the infinite-regular tree 
$$\u=\bigsqcup_{n\geq0}\N^n,$$
where $\N^n$ is the set of words with $n$ letters, and by convention
$\N^0=\{\varnothing\}$. For two words $u,v$, we let $uv$ be their
concatenation and $|u|,|v|$ their length (with the convention
$|\varnothing|=0$). If $u$ is a word and $A\subseteq \u$ we let
$uA=\{uv:v\in A\}$, and say that $u$ is a prefix of $v$ if $v\in u\u$,
in which case we write $u\vdash v$.  A {\em planar tree} is a finite
subset $\bt$ of $\u$ such that
\begin{itemize}
\item $\varnothing\in \bt$, it is called the {\em root} of $\bt$,
\item for every $u\in \u$ and $i\in\N$, if $ui\in \bt$ then $u\in
  \bt$, and $uj\in\bt$ for every $1\leq j\leq i$.
\end{itemize}
We let $\TT$ be the set of planar trees, which we simply refer to as
{\em trees} in the sequel.  For a tree $\bt\in \TT$ and $u\in\bt$, the
number $c_\bt(u)=\max\{i\in\Z_+:ui\in\bt\}$, with the convention
$u0=u$, is the number of children of $u$. We
say that $v$ is an ancestor of $u$ if $v\vdash u$.  An element
$u\in\bt$ is called a {\em vertex} of $\bt$, and the length $|u|$ of
the word $u$ is called the {\em height} of $u$ in $\bt$. The vertices
of $\bt$ with no children are called {\em leaves}.  For $\bt$ a planar
tree and $u$ a vertex of $\bt$, we let $\bt_u=\{v\in {\cal
  U}:uv\in\bt\}$, and we call it the {\em fringe subtree} rooted at
$u$ (it is trivially checked that it is indeed a tree). The `remaining
part' $[\bt]_u= \{u\}\cup(\bt\setminus u\bt_u)$ is called the {\em
 subtree of $\bt$ pruned at $u$}.
Any planar tree $\bt$ is endowed with the linear order which is the
restriction to $\bt$ of the usual lexicographical order $\prec$ on
$\u$ ($u\prec v$ if $u\vdash v$ or if $u=wu',v=wv'$ where
$u'_1<v'_1$). We call it the {\em depth-first order}.

In addition to trees, we will consider forests, which are defined as
nonempty subsets of $\u$ of the form
$$\ff=\bigcup_k k\bt_{(k)},$$ 
where $(\bt_{(k)})$ is a finite or infinite sequence of trees, which
are called the {\em components} of $\ff$. We let $\FF$ be the set of
forests.  The quantity $c_\ff(u)$ (number of children of $u\in\ff$) is
defined as for trees, and we let
$c_\ff(\varnothing)\in\N\sqcup\{\infty\}$ be the number of tree
components of $\ff$. We define $\ff_u=\{v:uv\in\ff\}\in\TT$ if
$u\in\ff$, and $\ff_u=\varnothing$ otherwise, so in particular
$\ff_k,1\leq k\leq c_\ff(\varnothing)$ are the tree components of
$\ff$. Also, let $[\ff]_u=\{u\}\cup\ff\setminus u\ff_u\in\FF$. If
$u\in \ff$ is a vertex of the forest, we call $|u|-1$ the height of
$u$. Notice that the notion of height of a vertex relies on whether we
are considering the vertex to belong to a tree or a forest, the reason
being that we want the roots $1,2,3,\ldots$ (or {\em floor}) of the forest
to be at height $0$. There should be no ambiguity according to the
context.

A {\em $K$-type planar tree}, or simply a multitype tree if the number
$K$ is clear from the context, is a pair $(\bt,e_\bt)$ where
$\bt\in\TT$ and $e_\bt:\bt\to[K]$.  We let $\TT^{(K)}$ be the set of
$K$-type trees. For $u\in\bt$, $e_\bt(u)$ is called the type of
$u$. If $(\bt,e_\bt)\in\TT^{(K)}$ and for $u\in\bt$, we let
$(\bt,e_\bt)_u\in\TT^{(K)}$ be the pair $(\bt_u,e_{\bt_u})$ where
$e_{\bt_u}(v)=e_\bt(uv)$. Similarly, $[(\bt,e_\bt)]_u$ is the tree
$[\bt]_u$ marked by the restriction of $e_\bt$ to $[\bt]_u$. Similar
definitions hold in a straightforward way for marked and $K$-type
forests $(\ff,e_\ff)$, whose set is denoted by $\FF^{(K)}$.

In the sequel, we will often denote the marking functions $e_\bt,e_\ff$
by $e$ when it is free of ambiguity, and will even denote elements of
$\TT^{(K)},\FF^{(K)}$ by $\bt$ or $\ff$, i.e.\ without explicitly
mentioning $e$. It will be understood then that
$\bt_u,\ff_u,...$ are marked with the appropriate function. We let,
for $i\in[K]$,
$$\TT^{(K)}_i=\{\bt\in\TT^{(K)}:e(\varnothing)=i\}\, ,$$ and for
$\bx=(x_j)$ a finite or infinite sequence with terms in $[K]$,
$$\FF^{(K)}_\bx=\{\ff\in \FF^{(K)}:e(j)=x_j \quad \forall j\}.$$ 
For $\bt\in\TT^{(K)}$ and $i\in[K]$, we let
$\bt^{(i)}=\{u\in\bt:e(u)=i\}$, and $\ff^{(i)}$ is the corresponding
notation for $\ff\in\FF^{(K)}$.

Let $\W_K=\bigsqcup_{n\geq 0}[K]^n$ be the set of finite, possibly
empty $[K]$-valued sequences, and consider the natural projection
$p:\W_K\to \Z_+^K$, where $p(\bw)=(p_i(\bw),i\in[K])$ and
$p_i(\bw)=\#\{j:\bw_j=i\}$ counts the number of elements of $\bw$
equal to $i$. Notice that for every $K$-type tree $\bt$, any $u\in
\bt$ determines a sequence $\bw_{\bt}(u)=(e(ui),1\leq i\leq
c_\bt(u))\in \W_K$ with length $|\bw_\bt(u)|=c_\bt(u)$. The vector
$p(\bw_\bt(u))$ counts the number of children of $u$ of each type.

Finally, we will frequently have to count the number of ancestors of
some vertex of a tree or a forest, that satisfy certain specific
properties. 
For $\bt\in\TT^{(K)}$ and $u\in\bt$, we let $\anc_\bt^u({\sf P}(v))$
be the number of ancestors $v\vdash u$ that satisfy the property ${\sf
  P}$.  For instance, $\anc_\bt^u(e(v)=i,c_\bt(v)=k)$ counts the
number of ancestors of $u$ with type $i$ and $k$ children.

\subsection{Galton-Watson trees}\label{sec:galton-watson-trees}

Let $\zzeta= (\zeta^{(i)},i\in[K])$ be a family of probability
measures on the set $\W_K$. We call $\zzeta$ an {\em ordered}
offspring distribution. It is said to be non-degenerate (resp.\
critical, resp.\ to have finite variance) if the family of measures
$\mmu=p_*\zzeta:=(p_*\zeta^{(i)},i\in[K])$ is non-degenerate (resp.\
critical, resp.\ satisfies (\ref{eq:6})), where $p_*\zeta^{(i)}$ is
the push-forward of $\zeta^{(i)}$ by $p$.

For $i\in[K]$, we now construct a distribution on $\TT^{(K)}_i$ such
that
\begin{itemize}
\item
different vertices have independent offspring, and 
\item type-$j$ vertices have a set of children with types given by a
  sequence $\bw\in\W_K$ with probability $\zeta^{(j)}(\bw)$.
\end{itemize}
To do this, let $({\bf W}_u^i=(W_u^i(l),1\leq l\leq |{\bf W}_u^i|),
1\leq i\leq K,u\in\u)$ be a $[K]\times\u$-indexed family of
independent random variables such that ${\bf W}_u^i$ has law
$\zeta^{(i)}$. Then recursively, construct a subset
$\mathfrak{t}\subset \u$ together with a mark $e:\mathfrak{t}\to[K]$
by letting $\varnothing\in\mathfrak{t},e(\varnothing)=i$, and if $u\in
\mathfrak{t},e(u)=j$, then $ul\in\mathfrak{t}$ if and only if $1\leq
l\leq |{\bf W}^j_u|$, and then $e(ul)=W^j_u(l)$.

It is straightforward to check that $\mathfrak{t}$ has the properties
of a planar tree, except that it might be infinite. Moreover, it is
straightforward from the construction that the process
$$\bZ_n(\mathfrak{t})=\left(\#\left\{u\in \mathfrak{t}:|u|=n,
    e(u)=i\right\},i\in[K]\right)\, , \qquad n\geq 0\, ,$$ is a
multitype GW process with offspring distribution
$\mmu=p_*\zzeta$, and started from a single type-$i$ individual. In
particular, under the criticality assumption $\varrho= 1$, this
process becomes extinct a.s., so that $\mathfrak{t}$ is finite a.s.\
and hence is a tree a.s.. In this case, we let $P^{(i)}_\mmu$, or
simply $P^{(i)}$, be the law of $\mathfrak{t}$ on $\TT^{(K)}_i$. The
probability measure $P^{(i)}$ is entirely characterized by the
formulae
$$P^{(i)}(T=\bt)=\prod_{u\in\bt}\zeta^{(e(u))}(\bw_{\bt}(u)),$$
where $T:\TT^{(K)}\to\TT^{(K)}$ is the identity map and $\bt$ ranges
over finite $K$-type trees.

Similarly, if $\bx=(x_j,1\leq j\leq r)\in \W_K$, we 
define $P^\bx$ as the image measure of $\bigotimes_{i=1}^r P^{(x_i)}$ 
by 
$$(\bt_{(1)},\ldots,\bt_{(k)})\longmapsto \bigcup_{1\leq k\leq r}k\bt_{(k)},$$
i.e., it is the law that makes the identity map
$F:\FF^{(K)}\to\FF^{(K)}$ the random forest whose tree components
$F_i,1\leq j\leq r$ are independent with respective laws
$P^{(x_i)}$. A similar definition holds for an infinite sequence
$\bx\in[K]^\N$. It will be convenient to use the notation $P^\bx_r$
for $P^{(x_1,\ldots,x_r)}$ when $\bx\in[K]^\N$.

\subsection{Convergence of height
  processes}\label{sec:main-result-1}

For $\bt\in\TT$, we let $\varnothing=u_\bt(0)\prec
u_\bt(1)\prec\ldots\prec u_\bt(\#\bt-1)$ be the list of vertices of
$\bt$ in depth-first order. When there is no ambiguity on $\bt$, we
simply denote them by $u(0),u(1),\ldots$. Let $(H_n^{\bt}=|u(n)|,n\geq
0)$ be the {\em height process} of $\bt$, with the convention that
$|u(n)|=0$ for $n\geq \#\bt$.  If $\ff\in\FF$, we similarly let
$u_\ff(0)\prec u_\ff(1)\prec\ldots$, or simply $u(0),u(1),\ldots$ be
the depth-first ordered list of its vertices ($u(0)=1$), and define
$(H_n^\ff,n\geq 0)$ by $H^\ff_n=(|u(n)|-1)\ind_{\{n\leq \#\ff-1\}}$
(again, the convention differs because we want the floor $1,2,\ldots$
of $\ff$ to be at height $0$).
Also, for $n\geq 0$ and $\ff\in\FF$,
let $\Upsilon^{\ff}_n$ be the first letter of $u(n)$
with the convention that for $n\geq\#\ff$, it equals the number of
components of $\ff$.

For $(\bt,e)\in\TT^{(K)}$ and $i\in[K]$, we let
$$\Lambda^\bt_i(n)=
\#\{0\leq k\leq n:e(u(k))=i\}$$ be the number of type-$i$ individuals
standing before the $n+1$-th individual in depth-first order. The
quantity $\Lambda^\ff_i$ is defined similarly for
$(\ff,e)\in\FF^{(K)}$.

\begin{thm}\label{T1}
  Let $\zzeta$ be an ordered offspring distribution such that
  $\mmu=p_*\zzeta$ satisfies {\rm (H)}. Recall the notations
  $\ba,\bb,\sigma$ around (\ref{sigma}). Then

\noindent{\rm (i)} 
Under $P^{\bx}$, for some arbitrary $\bx\in [K]^{\N}$, the following
convergence in distribution holds for the Skorokhod topology on the
space $\D(\R_+,\R)$ of right-continuous functions with left limits:
$$\left(\frac{H^{F}_{\lfloor ns\rfloor}}{\sqrt{n}},s\geq 0\right)
\build\longrightarrow_{n\to\infty}^{d} \left(\frac{2}{\sigma}|B_s|
  ,s\geq 0\right)\, ,$$ where $B$ is a standard one-dimensional
Brownian motion.

\noindent{\rm (ii)} For every $i\in[K]$, if ${\bf i}$ is the constant
sequence $(i,i,\ldots)$, then, under $P^{{\bf i}}$, the following
convergence in distribution in $\D(\R_+,\R)$ holds jointly with that
of {\rm (i)}:
$$\left( \frac{\Upsilon^{F}_{\lfloor ns\rfloor}}{\sqrt{n}}, s\geq 0\right)
\build\longrightarrow_{n\to\infty}^d \left( \frac{\sigma}{b_i}L^0_s,
  s\geq 0\right),$$ where $(L^0_t,t\geq 0)$ is the local time of $B$
at level $0$, normalized as the density of the occupation measure of
$B$ at $0$ before time $t$.

\noindent{\rm (iii)} Moreover, for any $\bx$, 
$$\left(\frac{\Lambda^F_i(\lfloor ns\rfloor)}{n},s\geq 0\right)
\build\longrightarrow_{n\to\infty}^{P^\bx}(a_i s,s\geq 0)\, 
,\qquad i\in[K]\, .$$ 
Here, the convergence is convergence in probability under $P^\bx$, for
the topology of uniform convergence over compact subsets of $\R_+$.
\end{thm}

Note that Theorem \ref{T1} could be also stated purely within a tree
formalism, without reference to height processes. In a rough way,
Theorem \ref{T1} says that multitype GW forests converge once properly
rescaled to a random forest, for a certain topology on the set of
tree-like metric spaces. This limiting forest is made of tree
components that are described by a Poisson process whose intensity
measure is the $\sigma$-finite Brownian continuum tree measure
$\Theta$ of \cite{duqlegprep}.

Let us comment on this result. First, (i) says that the height process
of a multitype forest, when properly rescaled, always looks the same
as a reflected Brownian motion with some prescribed scale factor,
whatever the roots are. Moreover, the scale factor depends only on
$\mmu=p_* \zzeta$, meaning that the exact way in which each set of
children is ordered does not affect the asymptotic distributional
shape of the forest. 

However, the value of $\bx$ actually matters if one wants to get a
closer look at the `limiting forest'.  Indeed, (ii) says that if one
wants to extract every single tree of a forest grown from a type-$i$
floor, one should proceed at a certain speed which does depend on
$i$. In particular, in a general mixed floor as in (i), when taking a
tree component in the limiting forest, one is in general unable to
recover the rank of the tree it comes from in the discrete picture.

Last, (iii) implies that
$$n^{-1}d\Lambda^F_i(\lfloor ns\rfloor)\ind_{\{s\geq 0\}}
\build\longrightarrow_{n\to\infty}^{P^\bx}a_i\, \d s\ind_{\{s\geq
  0\}}\, , \qquad 1\leq i\leq K\, ,$$ for the topology of vague
convergence of measures, hence showing that $\ba$ provides the
asymptotic relative weights of different types, which are not
influenced by the types of the roots. This is known as the {\em
  convergence of types theorem}, see \cite{AthNey}, and it is reproven
by different methods in Proposition \ref{J} below. Theorem \ref{T1}
(iii) gives the extra information that all types are homogeneously
distributed in the limiting
tree. 

We mention that the topology for the weak convergence of (i) and (ii)
could simply be the the uniform topology over compact subsets, since
all limits are continuous.

We also obtain as a corollary the following theorem of \cite{VD}, in
the more general case of irreducible mean matrix (\cite{VD} is in the
aperiodic case). 

\begin{crl}\label{sec:main-result}
  Let $\HT(\bt)$ be the maximal height of a vertex in $\bt$. Under the
  same assumptions, as $n\to\infty$, we have
$$nP^{(i)}(\HT(T)\geq n)\build\longrightarrow_{n\to\infty}^{}
\frac{2b_i}{\ba\cdot\bQ(\bb)}.$$
\end{crl}

Conditioned versions of Theorem \ref{T1} also hold. We say that $\mmu$
(or $\zzeta$) has small exponential moments if there exists $\eps>0$
such that
\begin{equation}\label{eq:7}
  \sup_{i\in[K]}
  \sum_{\bz\in \Z_+^K}\exp(\eps|\bz|_1)\mu^{(i)}(\bz)
=\sup_{i\in[K]}
  \sum_{\bw\in \W_K}\exp(\eps|\bw|)\zeta^{(i)}(\bw)<\infty\, ,
\end{equation}
where $|\bz|_1=z_1+\ldots+z_K$.  In the following statement, as well
as in all statements in the paper involving conditioned laws, we make
the assumption that $n$ goes to infinity along some subsequence, so
that all the conditioning events that are considered have positive
probabilities.

\begin{thm}\label{sec:convergence-gw-trees}
  Assume that hypothesis {\rm (H)} holds and that $\mmu$ has small
  exponential moments. Then for every $i,j\in[K]$, the following
  convergence in distribution holds on $\D([0,1],\R)$:
$$
\left(\frac{H^T_{\lfloor \#T t\rfloor}}{n^{1/2}},0\leq t\leq 1\right)\quad \mbox{
  under }\quad P^{(i)}(\cdot|\#T^{(j)}=n)
\build\longrightarrow_{n\to\infty}^{d}
\left(\frac{2}{\sigma\sqrt{a_j}}B^{{\rm ex}}_t,0\leq t\leq 1\right)\,
,$$ where $B^{{\rm ex}}$ is the standard Brownian excursion with
duration $1$. 

Moreover, $((\#T)^{-1}\Lambda^T_k(\lfloor \#T t\rfloor),0\leq t\leq 1)$ under
$P^{(i)}(\cdot|\#T^{(j)}=n)$ converges in probability to $(a_k t,0\leq
t\leq 1)$ for the uniform norm as $n\to\infty$, for every
$i,j,k\in[K]$.
\end{thm}

The main ingredient used to prove Theorems \ref{T1} and
\ref{sec:convergence-gw-trees} is the following result:

\begin{prp}[Corollary 2.5.1 in \cite{duqleg02}, Theorem 3.1 in
  \cite{duq02}]\label{thduqleg}  
  Theorems \ref{T1} and \ref{sec:convergence-gw-trees} are true in the
  case $K=1$.
\end{prp}

To be completely accurate, Theorem 2.5.1 in \cite{duqleg02} gives the
convergence in distribution of $n^{-1/2}(V^F_{\lfloor
  ns\rfloor},H^F_{\lfloor ns\rfloor},s\geq 0)$ to $(\sigma
B_s,2(B_s-I_s)/\sigma,s\geq 0)$, where $B$ is a standard Brownian
motion with infimum process $I_s=\inf_{0\leq u\leq s}B_u$, and where
$V^F$ is the so-called {\L}ukaciewicz walk of $F$ (see the proof of
Proposition \ref{sec:proof-theorem-reft2} for the definition), whose
only property we need at this point is that $\inf_{0\leq k\leq
  n}V^F_k=-\Upsilon^F_n$. This entails that
$n^{-1/2}(\Upsilon^F_{\lfloor ns\rfloor},H^F_{\lfloor ns\rfloor},s\geq
0)$ converges in distribution to $(- \sigma
I_s,2(B_s-I_s)/\sigma,s\geq 0)$, which by L\'evy's theorem has same
distribution as $(\sigma L^0_s,2|B_s|/\sigma,s\geq 0)$.

Also, Duquesne \cite{duq02} makes the assumption that the offspring
distribution is aperiodic to avoid conditioning events of zero
probability, but the proofs still work by considering subsequences as
we do.  The conditioned version of Proposition \ref{thduqleg} was
first stated in Aldous \cite{aldouscrt93}, using the so-called {\em
  contour} process, rather than the height process, to encode the
discrete trees. We stress that in \cite{aldouscrt93,duq02}, the
authors do not assume that the offspring distribution has small
exponential moments, and we expect Theorem
\ref{sec:convergence-gw-trees} to hold without this extra
hypothesis. As a matter of fact, by making occasional changes in the
proofs below, one can show that a sixth moment for $\mmu$ is
sufficient, and we suspect that a second moment condition is enough.

We finally stress that Proposition \ref{thduqleg} is in fact a
particular case of Duquesne and Le Gall's results
\cite{duqleg02,duq02}, which deal with the case of offspring
distributions belonging to other stable domains of attraction than the
Gaussian one.

The idea of the proof of Theorems \ref{T1},
\ref{sec:convergence-gw-trees} will be to use an inductive argument on
$K$ in order to apply the foregoing proposition. We rely strongly on
an ancestral decomposition for multitype trees adapted from
\cite{KLPP} and to be developed in Section \ref{sec:ancestral}, in
which we also state some facts about monotype trees and the
Perron-Frobenius theorem. The proof of Theorem \ref{T1} is then given
in Sections \ref{sec:convtypes} and \ref{sec:proof-i-ii}.

\subsection{Convergence of multitype snakes to the Brownian
  snake}\label{snakes}

Let us now couple the multitype branching process with a spatial
motion. As in \cite{jfmgm05}, we are interested in the case where the
motion, which is parametrized by the vertices of the tree, has a step
distribution around a given vertex that may depend locally on the
tree, through the type and children of the vertex. We stress that our
method could most likely be applied to other kinds of step
distributions that depend on the structures of neighborhoods of the
vertices which are not too large (i.e.\ have uniformly negligible
scale compared to the height of the large trees).

Consider a family of probability distributions $\nu_{i,\bw}$
respectively on $\R^{|\bw|}$, and indexed by types $i\in[K]$ and
$\bw\in\W_K$.  For $(\bt,e)$ a multitype tree and every $u\in \bt$
with type $e(u)=i$ and children vector $\bw_\bt(u)=\bw$, we take a
random variable $(Y_{uj},1\leq j\leq |\bw|)$ with law $\nu_{i,\bw}$,
independently over distinct vertices. Let $\nu_\bt$ be the law of the
random vector $(Y_u,u\in\bt)$ thus obtained, with the convention that
$Y_\varnothing=0$. We let
$$\T^{(K)}=\{(\bt,e,(y_u,u\in\bt)):(\bt,e)\in\TT^{(K)},{\bf y}\in \R^{\bt}\}$$
be the set of multitype trees with {\em spatial marks} on the
vertices, and let $\P^{(i)}$ be the probability measure $\d
P^{(i)}(\bt,e)\otimes \d\nu_\bt( {\bf y})$.

Similarly, for $(\ff,e)$ a multitype forest we let $\nu_\ff$ be the
law of the random variable $(Y_u,u\in\ff)$, where $Y_n=0$ for any
$n\in\N$ and the random vectors $(Y_{ui}:1\leq i\leq c_\ff(u))$
for $u\in \ff$ are independent with respective laws
$\nu_{e(u),\bw_\ff(u)}$. For $\bx=(x_1,x_2,\ldots)$ a finite or
infinite sequence of types we let $\P^\bx=\d P^{\bx}((\ff,e))\otimes
\d\nu_\ff({\bf y})$, which is a probability distribution on the set
$$\F^{(K)}=
\{(\bt,e,(y_u,u\in\bt)):(\bt,e)\in\TT^{(K)},{\bf y}\in \R^{\bt}\}.$$

Notice that in the definitions of $\P^{(i)}$ and $\P^{\bx}$, only the
measures $\nu_{i,\bw}$ for $(i,\bw)$ such that $\zeta^{(i)}(\bw)>0$
matter. By convention, we let $\nu_{i,\bw}$ be the Dirac mass at
$0\in\R^{|\bw|}$ for all irrelevant indices. We say that the family
$(\nu_{i,\bw},i\in[K],\bw\in\W_K)$ is non-degenerate if $\nu_{i,\bw}$
is not a Dirac mass for at least one of the relevant indices
$(i,\bw)$, so that there is `some randomness' in the spatial
displacement. We say that $\nnu=(\nu_{i,\bw},i\in[K],\bw\in\W_K)$ is
centered if all distributions $\nu_{i,\bw}$ are. 

For $(\bt,e,{\bf y})\in\T^{(K)}$, define
$$S^{\bt,e,{\bf y}}_u=\sum_{v\vdash u}y_{v},$$
and let $S^{\bt}_{k}$ stand (a little improperly, but for lighter
notations) for $S^{\bt,e,{\bf y}}_{u(k)}$, with the convention that it
equals $0$ for $k\geq \#\bt$.
A similar definition holds for $S^\ff$ where $(\ff,e,{\bf y})\in
\F^{(K)}$, where we use the convention $y_\varnothing=0$. In the
following statement and in the sequel, $|\cdot|_2$ is the Euclidean
norm of the vector $x$.

\begin{thm}\label{T2}
  Assume $\mmu=p_*\zzeta$ satisfies {\rm (H)} and admits some
  exponential moments, and that $\nnu$ is non-degenerate and centered.
  Suppose also that every $\nu_{i,\bw}$ admits a moment ${\tt
    M}_{i,\bw}=\langle\nu_{i,\bw},|y|_2^{8+\xi}\rangle$, for some
  $\xi>0$, such that
  \begin{equation}\label{eq:5}
    \sup_{i\in[K]}{\tt M}_{i,\bw}=O(|\bw|^D)
\end{equation}
for some $D>0$ (here $|y|_2$ stands for the Euclidean norm of
$y$). Write
$$\Sigma=\sqrt{\sum_{i\in[K]}a_i\sum_{\bw\in\W_K}\zeta^{(i)}(\bw)
  \sum_{j=1}^{|\bw|}b_{w_j}\langle\nu_{i,\bw},y_j^2\rangle}\,
\in(0,\infty)\, .$$ Then for any $\bx\in[K]^\N$, under $\P^{\bx}$, the
following convergence in distribution on $\D(\R_+,\R)$ holds jointly
with that of {\rm (i)} in Theorem \ref{T1}:
$$\left(
  \frac{1}{n^{1/4}}S^F_{\lfloor ns\rfloor},s\geq
  0\right)\build\longrightarrow_{n\to\infty}^{d}\left(
  \Sigma\sqrt{\frac{2}{\sigma}}R_s,s\geq 0\right),$$ where
conditionally on the Brownian motion $B$ of {\rm (i)} in Theorem
\ref{T1}, $R$ is a Gaussian process with covariance
$$\cov(R_s,R_{s'})=\inf_{s\wedge s'\leq u\leq s\vee s'}
|B_u|.$$ If $\bx=(i,i,\ldots)$ for some $i\in[K]$, then the
convergence holds jointly with that of {\rm (ii)} in Theorem
\ref{T1}.
\end{thm}

The analogous statement for conditioned laws is: 

\begin{thm}\label{sec:conv-mult-snak}
Under the same hypotheses as Theorem \ref{T2},  
the following convergence in distribution on
$\D([0,1],\R)$ holds jointly with that of Theorem
\ref{sec:convergence-gw-trees}: for every $i,j\in[K]$,
$$
\left(\frac{S^T_{\lfloor \#T t\rfloor}}{n^{1/4}},0\leq t\leq 1\right)\quad \mbox{
  under }\quad \P^{(i)}(\cdot|\#T^{(j)}=n)
\build\longrightarrow_{n\to\infty}^{d}
\left(\Sigma\sqrt{\frac{2}{\sigma\sqrt{a_j}}}R^{{\rm ex}}_t,0\leq
  t\leq 1\right)\, ,$$ where conditionally on the Brownian excursion
$B^{{\rm ex}}$ of Theorem \ref{sec:convergence-gw-trees}, $R^{{\rm
    ex}}$ is a Gaussian process with covariance
$$\cov(R^{{\rm ex}}_s,R^{{\rm ex}}_{s'})=\inf_{s\wedge s'\leq u\leq s\vee s'}
B^{{\rm ex}}_u\, .$$
\end{thm}

\rem
By contrast with \cite{janmarck05} and \cite{jfmgm05}, our
hypothesis on the spatial displacement is an $8+\xi$-moment assumption
rather than a $4+\xi$-moment assumption. We believe such a weaker
hypothesis to be sufficient, but were not able to prove it,
essentially because we could not prove what we believe to be the best
H\"older norm bounds in Proposition \ref{sec:proof-theorem-reft2}. See
the remark after the latter's statement.

\bigskip

\noindent{\bf Note. }In this paper, we will make a frequent use of
exponential bounds for real sequences $(x_n,n\geq 0)$, namely that
$|x_n|\leq \exp(-n^\eps)$ for some $\eps>0$ and large enough $n$. To
simplify notations and avoid referring to changing $\eps$'s, we write
$x_n=\oe(n)$ in this case.

\section{Proof of Theorem \ref{T1}}\label{sec:prelimilary-results}

\subsection[Ancestral decompositions]{Ancestral decomposition for 
multitype Galton-Watson trees and forests}\label{sec:ancestral}

Let $\zzeta=(\zeta^{(1)},\ldots,\zeta^{(K)})$ be a non-degenerate
critical ordered offspring distribution. For $i\in[K]$, define the
size-biased measure
$$\wh{\zeta}^{(i)}(\bw)=\frac{p(\bw)\cdot\bb}{b_i}\zeta^{(i)}(\bw)\, ,\qquad 
\bw\in\W_K\, ,$$ and notice that these are probability measures on
$\W_K$ by the definition of $\bb$, since they are pushed by $p$ to the
measure
$$\wh{\mu}^{(i)}(\bz)=\frac{\bz\cdot \bb}{b_i}\mu^{(i)}(\bz)\, ,\qquad
\bz\in\Z_+^K\, ,$$ and they do not charge the null sequence
$\varnothing$.  On some probability space $(\Omega,\a,P)$, let
$((\bW^i_u,\wh{\bW}^i_u,l^i_u,i\in[K]),u\in\u)$ be a family of
$\u$-indexed independent random vectors, such that $\bW^i_u$ has law
$\zeta^{(i)}$, $\wh{\bW}^i_u$ has law $\wh{\zeta}^{(i)}$, and
$$P(l^i_u=k|\wh{\bW}^i_u)=\frac{b_{\wh{W}^i_u(k)}}{p(\wh{\bW}^i_u)\cdot
  \bb}\, .$$ Otherwise said, $l^i_u$ is equal to $k$ with probability
proportional to $b_{w_k}$ given $\wh{\bW}^i_u=\bw$.  Fix some
$i\in[K]$. Recursively, we build a set $\wh{\mathfrak{t}}$, a mark
$\wh{e}:\mathfrak{t}\to[K]$ and a sequence $(V_n,n\geq 0)$ by first
letting $V_0=\varnothing\in\wh{\mathfrak{t}}$,
$\wh{e}(\varnothing)=i$, and given
$V_0,\wh{e}(V_0),\ldots,V_n,\wh{e}(V_n)$ have been constructed with
$\wh{e}(V_n)=j$, we let $V_{n+1}=V_nl^j_{V_n}$ and
$\wh{e}(V_{n+1})=\wh{W}^j_{V_n}(l^j_{V_n})$. Then, for
$u\in\wh{\mathfrak{t}}$ with $\wh{e}(u)=j$,
\begin{itemize}
\item if $u=V_n$ for some $n$, then $ul\in\wh{\mathfrak{t}}$ if and only if
  $1\leq l\leq |\wh{\bW}^j_u|$. For such $l$,
  $\wh{e}(ul)=\wh{W}^j_u(l)$, and
\item otherwise, $ul\in\wh{\mathfrak{t}}$ if and only if $1\leq l\leq
  |\bW^j_u|$. For such $l$, $\wh{e}(ul)=W^j_u(l)$.
\end{itemize}
The set $\wh{\mathfrak{t}}$ thus obtained has the properties of a
tree, except that it is infinite. More precisely, it consists of an
infinite `spine' $V_0,V_1,\ldots$ of distinguished vertices, which is
interpreted as an infinite ancestral line along which individuals of
type $i$ have a $\wh{\zeta}^{(i)}$-distributed offspring sequence,
among which each element $j$ is selected as a distinguished successor
with probability proportional to $b_j$. Then, non-distinguished
individuals have a regular GW descendence with offspring distribution
$\zzeta$. It is easy to prove that the spine is the unique infinite
simple path of vertices in $\wh{\mathfrak{t}}$ starting from
$\varnothing$.

We let $\wh{P}^{(i),h}$ be the law of $([\wh{\mathfrak{t}}]_{V_h},V_h)$.  It is
a distribution on the set $\{(\bt,e,v):(\bt,e)\in\TT^{(K)},v\in \bt\}$
of pointed trees, on which we let $(T,V)$ be the identity map.  For
any finite sequence $\bx$ of types, let also $\wh{P}^{\bx,j,h}$ be the
law under which $(F_l,l\neq j), F_j$ are pairwise independent with
respective laws $(P^{(x_l)},l\neq j), \wh{P}^{(x_j),h}$.

Notice that by construction, the types of the distinguished
individuals $V_0,V_1,\ldots$ form a Markov chain,
whose transition law is
\begin{equation}\label{transition}
  \wh{P}^{(i)}(e(V_{h+1})=j'|e(V_h)=j)=\sum_{\bz\in \Z_+^K}\frac{
    \wh{\mu}^{(j)}(\bz)b_{j'}z_{j'}}{\bb\cdot\bz}=\frac{b_{j'}}{b_j}m_{jj'}.
\end{equation}
The stationary distribution of this Markov chain is easily checked to
be the vector $(a_ib_i,1\leq i\leq K)$, because of the normalization
$\ba\cdot\bb=1$.  This will be useful in the sequel.

\begin{lmm}\label{ancestral}
  For any finite sequence of types $\bx=(x_j,1\leq j\leq r)$, and any
  non-negative functions $G_1,G_2$,
\begin{equation}\label{ancestreq}
  E^{\bx}\left[\sum_{v\in F}G_1(v,[F]_v)G_2(F_v)\right]
  =\sum_{j=1}^r b_{x_j} \sum_{h\geq 0}\wh{E}^{\bx,j,h}\left[\frac{G_1(V,F)
      E^{(e(V))}\left[G_2(T)\right]}{b_{e(V)}}\right]. 
\end{equation}
\end{lmm}

\proof Let $(\ff,e_\ff)$ be a $K$-type forest, $u$ a leaf of $\ff$
(i.e.\ a vertex with no child) and $(\bt,e_\bt)$ a $K$-type tree with
$e_\bt(\varnothing)=e_\ff(u)$. Then it is enough to show the result
for $G_1(u,\ff)G_2(\bt)=\ind_{\{u,\ff,\bt\}}$, as one can then use
linearity of the expectation.  In this case, the left-hand side of
(\ref{ancestreq}) is equal to $P^{\bx}(F=[\ff,u,\bt])$, where
$[\ff,u,\bt]$ is the only forest $(\ff',e')\in \FF^{(K)}$ containing
$u$ with $[\ff']_u=\ff$ and $\ff'_u=\bt$. We let $S=\{v:v\vdash
u,v\neq u\}$ and $A=\{v:v\notin S\cup\{ u\}, \neg v\in S\}$, where
$\neg v$ is the father of $v$, i.e.\ the word $v$ with its last letter
removed.  We can redisplay
$$P^{\bx}(F=[\ff,u,\bt])=\prod_{v\in \ff'}\zeta^{(e_\ff(v))}(\bw_{\ff'}(v)),$$
as
$$\prod_{v\in \bt}\zeta^{(e_\bt(v))}(\bw_\bt(v))
\prod_{v\in\ff,v\nvdash u}\zeta^{(e_\ff(v))}(\bw_\ff(v)) \prod_{v\in
  S}\zeta^{(e_\ff(v))}(\bw_\ff(v)).$$ In this expression, one can
factorize out the product of probabilities of the subtrees
$\ff'_u=\bt$ and $\ff'_j=\ff_j$ for $j\in\{1,\ldots,r\}\setminus
\{u_1\}$, $u_1$ the first letter of $u$, and $\ff_{v}$ for $v\in
A$. Letting $j=u_1$ and $u=ju'$, this shows
\begin{eqnarray*}
  \lefteqn{P^{\bx}(F=[\ff,u,\bt])}\\
  &=& P^{(e_\ff(u))}(T=\bt)
  \prod_{l\neq j,1\leq l\leq r}P^{(x_l)}(T=\ff_l)
  \prod_{v\in A}P^{(e_\ff(v))}(T=\ff_v)
  \prod_{v\in S}\zeta^{(e_\ff(v))}(\bw_\ff(v)).
\end{eqnarray*}
We can also rewrite the last product as
$$\prod_{v\in S}\wh{\zeta}^{(e_\ff(v))}(\bw_\ff(v))\frac{b_{e_\ff(v)}}{
  p(\bw_\ff(v))\cdot\bb},$$ and we finally recognize
\begin{eqnarray*}
  \lefteqn{ P^{\bx}(F=[\ff,u,\bt])}\\
  &=&\frac{b_{x_j}}{b_{e_\ff(u)}}\wh{P}^{(x_j),|u|-1}((T,V)=(\ff_j,u'))
  P^{(e_\ff(u))}(T=\bt)
  \left(\prod_{l\neq j,1\leq l\leq r}P^{(x_l)}(T=\ff_l)\right)  \\
  &=&b_{x_j} 
  \wh{E}^{\bx,j,|u|}[\ind_{\{u,\ff\}}(V,F)E^{(e(V))}[\ind_{\{\bt\}}(T)]/
  b_{e(V)}],
\end{eqnarray*}
which yields the result. \cq

\subsection{Around Perron-Frobenius'
  Theorem}\label{sec:around-perr-frob}

Let us first recall the well-known Perron-Frobenius Theorem, which in
this form can be found in \cite{seneta81}.
\begin{prp}[Perron-Frobenius]\label{frobenius}
  Let $M\in\mathfrak{M}_K(\R_+)$ be an irreducible matrix.

\noindent{\rm (i)} 
The matrix $M$ has a real eigenvalue $\varrho$ with maximal modulus,
which is positive, simple (i.e.\ it is a simple root of the
characteristic polyn\^omial of $M$), and every $\varrho$-eigenvector
has only non-zero entries, all of the same sign.

\noindent{\rm (ii)} Any eigenvector of $M$ with non-negative entries
is a $\varrho$-eigenvector, and hence has only positive entries.
\end{prp}

From this, we deduce the following useful

\begin{lmm}\label{spmii}
{\rm (i)} Suppose $M\in \mathfrak{M}_K(\R_+)$ is irreducible. 
Then its spectral radius $\varrho$ satisfies $\varrho>m_{ii}$, for all
$1\leq i\leq K$.

\noindent{\rm (ii)}
Suppose $K>1$ and $\varrho=1$.  Then the matrix
$\wt{M}\in\mathfrak{M}_{K-1}(\R_+)$ with entries
$$\wt{m}_{ij}=m_{ij}+\frac{m_{iK}m_{Kj}}{1-m_{KK}}\, ,\qquad 1\leq i,j\leq 
K-1$$
is also irreducible with spectral radius $1$.
\end{lmm}

\proof (i) is immediate by writing $\sum_{j=1}^Km_{ij}u_j=\varrho
u_i,1\leq i\leq K$, where ${\bf u}$ is a right $\varrho$-eigenvector
of $M$ with positive entries.

(ii) First, notice that the irreducibility of a matrix with
nonnegative entries only depends on which entries are non-zero, and
not on what their actual value is. Since $1-m_{KK}>0$ by (i), it is
sufficient to show that the matrix with entries
\begin{equation}\label{eq:21}
\max(m_{ij},m_{iK}m_{Kj})\, ,\qquad 1\leq i,j\leq K-1
\end{equation}
is irreducible when $M$ is. Assuming $M$ irreducible and upon
replacing $m_{ij}$ by $\ind_{\{m_{ij}>0\}}$, we may assume $M$ is the
adjacency matrix of a connected non-oriented graph on the vertices
$1,\ldots,K$. But now, the matrix of (\ref{eq:21}) is the adjacency
matrix of the graph on the vertices $1,\ldots, K-1$, where $i$ is
adjacent to $j$ if and only if they are either adjacent or both
adjacent to $K$ in the initial graph. It is straightforward to see
that this graph remains connected, so its adjacency matrix is
irreducible.

It remains to show that $\wt{M}$ has spectral radius $1$. First, it is
immediate that if $\ba,\bb$ are left and right $1$-eigenvectors of
$M$, then $(a_1,\ldots,a_{K-1})$ and $(b_1,\ldots,b_{K-1})$ are left
and right $1$-eigenvectors of $\wt{M}$. Indeed, for $1\leq j\leq K-1$,
$$\sum_{i=1}^{K-1}a_i\wt{m}_{ij}=a_j-a_Km_{Kj}+a_K(1-m_{KK})
\frac{m_{Kj}}{1-m_{KK}}=a_j, \mbox{ and }
\sum_{j=1}^{K-1}b_j\wt{m}_{ij}=b_i$$ for $1\leq i\leq K-1$. This is
enough to conclude by (ii) in Proposition \ref{frobenius}. \cq

\subsection{Reduction of trees}\label{sec:reduction-trees}

\subsubsection{A projection on monotype trees}\label{sec:proj-monotype-trees}

We describe a projection function $\Pi^{(i)}$ that goes from the set
of $K$-type planar forests to the set of monotype planar forests, and
which intuitively squeezes generations, keeping only the type-$i$
individuals.

Precisely, if $(\ff,e)$ is a $K$-type forest, we first let $v_{1}\prec
v_{2}\prec\ldots$ be the ordered list of vertices of $\ff^{(i)}$ such
that all ancestors of $v_k$ have type different from $i$. We consider
a forest $\Pi^{(i)}(\ff)=\ff'$ with as many tree components as there
are elements in $\{v_1,v_2,\ldots\}$. Thus, we start with the set of
roots $1,2,\ldots$ of $\ff'$. Then recursively, for each $u\in \ff'$,
let $v_{u1},\ldots,v_{uk}$ be the vertices of
$(v_u\ff_{v_u})\setminus\{v_u\}$ such that
\begin{itemize}
\item
$e(v_{uj})=i, 1\leq j\leq k$,
\item for every $1\leq j\leq k$, if $v_{uj}=v_ul_1,\ldots l_h$, then
  $e(ul_1\ldots l_r)\neq i$ for all $1\leq r< h$, and
\item
$v_{u1},\ldots,v_{uk}$ are arranged in lexicographical order.
\end{itemize}
Then, we add the vertices $u1,\ldots,uk$ to $\ff'$, and
continue iteratively. If $u\in\ff'$ has children $u1,\ldots,uk$, we
let 
$$N^{(i)}(u)=\#\left(\ff_{v_u}\setminus 
  \left(\bigcup_{j=1}^k\ff_{v_{uj}}\right)\right)-1\, ,$$ be
the number of vertices that have been deleted between $u$ and its children
during the operation. We also let 
$$N^{(i)}_{\rm floor}(n)=\#\{v\in \ff_n:e(w)\neq i\, \forall w\vdash
v\}\, ,$$ be the number of vertices of the $n$-th tree component of
$\ff$ that lie below the first layer of type-$i$ vertices.

If $(\bt,e)\in\TT^{(K)}$, we may apply the map $\Pi^{(i)}$ to the
forest $1\bt$, and get as a result a forest $\ff'$. We denote this
forest by $\Pi^{(i)}(\bt)$ with a slight abuse of notation. Notice
that $\Pi^{(i)}(\bt)$ has one tree component if $e(\varnothing)=i$, in
which particular case we denote it by $\ff'_1=\Pi^{(i)}(\bt)$, with
further abuse of notations. 

If $\mmu$ is an unordered offspring distribution and $i,j\in[K]$, we
let $\ov{\mu}_j^{(i)}$ be the distribution of
$c_{\Pi^{(i)}(T)}(\varnothing)$ under $P^{(j)}$ with the above
conventions: if $i\neq j$ then this counts the number of components of
the reduced forest, while if $i=j$ this is the number of children of
the root of the reduced tree. 

\begin{prp}\label{monotype}
{\rm (i)}
Let $\bx\in[K]^\N$. Under the law $P^\bx$, the forest $\Pi^{(i)}(F)$
is a (monotype) GW forest with offspring distribution
$\ov{\mu}_i^{(i)}$,  
whose mean and variance are equal to
$$\sum_{z\geq 0}z\ov{\mu}^{(i)}_i(\{z\})=1\qquad\mbox{and}\qquad
\var(\ov{\mu}_i^{(i)})=\frac{\sigma^2}{a_ib_i^2}.$$ In particular, it
is critical with finite variance.

\noindent{\rm (ii)}
For each $i$ and still under $E^\bx$, the random variables
$(N^{(i)}(u(n)),n\geq 0),(N^{(i)}_{\rm floor}(n),n\geq 1)$ are all
independent, and the variables $(N^{(i)}(u(n)),n\geq 0)$ are i.i.d.
The number of vertices which are deleted in the operation $\Pi^{(i)}$
between two generations has mean $$E^{\bx}\left[N^{(i)}(u(0))
\right]= \frac{1}{a_i}-1\, ,$$ and finite variance. Similarly, the
random variables $(N^{(i)}_{\rm floor}(n),n\geq 1)$ have finite
($x_n$-dependent) variance, as well as the laws $\ov{\mu}_i^{(j)}$ for
$i,j\in[K]$.   

\noindent{\rm (iii)} More generally, if $\mmu$ admits a finite $p$-th
moment with $p\in\N$ (resp.\ admits some exponential moments), then so
do $\ov{\mu}^{(j)}_i$ and the variables $N^{(i)}(u(n)),N^{(i)}_{\rm
  floor}(n+1),n\geq 0$ under $P^\bx$.
\end{prp}

The GW property of $\Pi^{(i)}(F)$ under $E^\bx$ is easy to obtain from
Jagers' theorem on stopping lines \cite{MR1014449}, each subtree
rooted at a vertex of type $i$ being a copy of the whole tree. This
also gives the independence statement in {\rm (ii)}. A detailed proof
of these intuitive statements would be cumbersome, so that we leave
the details to the interested reader, whom we refer to
\cite{MR1014449}.

The rest of the proof of this proposition will be done by removing
types one by one, and using an induction argument.

\subsubsection{From $K$ to $K-1$ types}\label{sec:from-k-k}

In this section, we will suppose that type $K$ is
deleted, keeping in mind that the general case is similar.  If
$\ff\in\FF^{(K)}$, we now let $v_1\prec v_2\prec\ldots$ be the ordered
list of vertices of $\ff$ such that $e(v_i)\neq K$ and $e(v)=K$ for
every $v\vdash v_i$. Recursively, given $v_u\in \ff$ has been
constructed, we let $v_{u1}\prec \ldots\prec v_{uk}$ be the
descendents of $v_u$ such that $e(v)=K$ for every $v_u\vdash v\vdash
v_{uj}$, $v\notin\{v_u,v_{uj}\}$, while $e(v_{uj})\neq K$. We then let
$\wt{\Pi}(\ff)$ be the set of $u\in{\cal U}$ such that $v_u$ has been
defined by our recursive construction, and naturally associate a type
$\wt{e}(u)=e(v_u)$ with $\wt{\Pi}(\ff)$. If $u\in\wt{\Pi}(\ff)$, we
let $\wt{N}(u)$ be the number of vertices of type $K$ that have been
deleted in this construction between $v_u$ and $v_{u1},\ldots,v_{uk}$,
namely,
$$\wt{N}(u)=\#\left(\ff_{v_u}\setminus\bigcup_{i=1}^k
  \ff_{v_{ui}}\right)-1\, ,$$ with the above notations. We also let 
$\wt{N}_{\rm floor}(n)$ be the number of vertices $v\in\ff^{(K)}_n$
that have only ancestors of type $K$.

\begin{lmm}\label{toutvabien}
  Let $\bx\in[K]^\N$. Then, under $P^\bx$:

{\rm (i)}
for any $1\leq i\leq K-1$,
  the forest $\wt{\Pi}(F)$ is a non-degenerate, irreducible, critical
  $K-1$-type GW forest. The (unordered) offspring
  distribution $\wt{\mmu}=(\wt{\mu}^{(i)},i\in[K-1])$ has generating
  functions
\begin{equation}\label{phij}
  \wt{\varphi}^{(i)}(\bs)=\varphi^{(i)}(\bs,\wt{\varphi}^{(K)}(\bs))\, ,
\end{equation}
for $1\leq i\leq K-1$ and $\bs\in[0,1]^{K-1}$, where $\wt{\varphi}^{(K)}$
is implicitly defined by 
\begin{equation}\label{phiK}
\wt{\varphi}^{(K)}(\bs)=\varphi^{(K)}(\bs,\wt{\varphi}^{(K)}(\bs))\, .
\end{equation}

{\rm (ii)} the $K-1$ sequences $(\wt{N}(u^{(i)}(n)),n\geq 0),i\in[K-1]$ are
independent and formed of i.i.d.\ elements, where $u^{(i)}(0)\prec
u^{(i)}(1)\prec\ldots$ is the ordered list of elements of $F$ with type
$i$. Their generating functions $\wt{\psi}^{(i)}$ respectively satisfy
\begin{equation}\label{eq:13}
  \wt{\psi}^{(i)}(s)=\varphi^{(i)}(1,\ldots,1,\wt{\psi}^{(K)}(s))\,
  ,\end{equation} 
where $\wt{\psi}^{(K)}$ is implicitly defined by
\begin{equation}\label{eq:14}
\wt{\psi}^{(K)}(s)=s\varphi^{(K)}(1,\ldots,1,\wt{\psi}^{(K)}(s))\, .
\end{equation}
The random variables $\wt{N}_{\rm floor}(n),n\geq 1$ are independent
as well.

{\rm (iii)} any integer or exponential moment conditions on $\mmu$ is
also satisfied by the laws $\wt{\mmu}$ and the random variables
$\wt{N}(u^{(i)}(n)),\wt{N}_{\rm floor}(n+1),n\geq 0,i\in[K-1]$. 
\end{lmm}

\proof {\rm (i)} Again, the GW property follows from the
construction. On some probability space, let $\wt{\bX}^{(j)}$ have
same distribution as the $\Z_+^{K-1}$-valued random vector of children
of a type-$j$ vertex of $\wt{\Pi}(F)$ under $P^\bx$. Then, by
separating the offspring of this vertex with types equal and different
from $K$, we obtain the identity in law
$$\wt{\bX}^{(j)}\build=_{}^d (X^{(j)}_1,\ldots,X^{(j)}_{K-1})
+\sum_{l=1}^{X^{(j)}_K}\wt{\bX}^{(K)}(l)\, ,$$ where $\bX^{(i)}$ has
distribution $\mu^{(i)}$ and is independent of
$(\wt{\bX}^{(K)}(l),l\geq 1)$, which are independent with same
distribution, and a vector $\wt{\bX}^{(K)}$ with this distribution
must satisfy
$$\wt{\bX}^{(K)}\build=_{}^d (X^{(K)}_1,\ldots,X^{(K)}_{K-1})
+\sum_{l=1}^{X^{(K)}_K}\wt{\bX}^{(K)}(l)\, ,$$ with similar
notations. These two expressions immediately translate as (\ref{phij})
and (\ref{phiK}).  Now let
$$\wt{m}_{ij}=\frac{\dd \wt{\varphi}^{(i)}}{\dd s_j}({\bf 1})\, ,\qquad 
1\leq i,j\leq K-1,$$ so that $\wt{M}=(\wt{m}_{ij})_{1\leq i,j\leq K}$
is the mean matrix associated with the $K-1$-type GW forest
$\wt{\Pi}(F)$ under $P^\bx$.  Differentiating (\ref{phij}) and
(\ref{phiK}) and letting $\bs$ increase to $(1,\ldots,1)$ gives
$$\wt{m}_{ij}=m_{ij}+\frac{m_{iK}m_{Kj}}{1-m_{KK}},$$
for $1\leq i,j\leq K-1$. Hence, $\wt{M}$ is defined as in Lemma
\ref{spmii}, so it is irreducible with spectral radius $1$.  Moreover,
the $K$-type GW process associated with $\wt{\Pi}(F)$ under $P^\bx$
has to be non-degenerate, because it dies in finite time a.s..

{\rm (ii)} The independence statement is again a consequence of
Jager's theorem, and Formulas (\ref{eq:13}) and (\ref{eq:14}) are
obtained by similar distributional equations arguments as above,
namely,
$$\wt{N}(u^{(i)}(n))\build=_{}^{d}\sum_{l=1}^{X^{(i)}_K}
\wt{N}^{(K)}_l\, ,$$ where $\bX^{(i)}$ has law $\mu^{(i)}$ and
$\wt{N}^{(K)}_l,l\geq 1$ are i.i.d.\ random elements independent of
$\bX^{(i)}$, that satisfy
$$\wt{N}^{(K)}_1\build=_{}^{d}1+ \sum_{l=1}^{X^{(K)}_K}
\wt{N}^{(K)}_l\, .$$ On the other hand, $\wt{N}_{\rm floor}(n)$ is
either $0$ or equal in distribution to $\wt{N}^{(K)}_1+1$ with the
same notation, according to whether $x_n\neq K$ or $x_n=K$.  

{\rm (iii)} is obtained by differentiating equations (\ref{phij}),
(\ref{phiK}), (\ref{eq:13}) and (\ref{eq:14}) $p$ times, while the
assertion on small exponential moments is obtained by applying the
implicit function theorem to the implicit functions
$\wt{\varphi}^{(K)},\wt{\psi}^{(K)}$. Details are left as an exercise
to the reader.  \cq

\bigskip

Notice that this provides an alternative way of showing that the
spectral radius of $\wt{M}$ is $\leq 1$, since the GW process has to
be (sub)-critical in order to become extinct a.s.  Recall that
$\ba,\bb$ are the left and right $1$-eigenvectors of $M$ with
$\ba\cdot {\bf 1}=1=\ba\cdot\bb$. In view of the proof of Lemma
\ref{spmii}, the left and right $1$-eigenvectors $\wt{\ba},\wt{\bb}$
of $\wt{M}$ satisfying $\wt{\ba}\cdot{\bf 1}=1=\wt{\ba}\cdot\wt{\bb}$
are given by
$$\wt{\ba}=\frac{1}{1-a_K}(a_1,\ldots,a_{K-1})\, ,\qquad 
\wt{\bb}=\frac{1-a_K}{1-a_Kb_K}(b_1,\ldots,b_{K-1}).$$ We are now
ready to give the

\bigskip

\noindent{\bf Proof of Proposition \ref{monotype}. }
(i) We prove this by induction on $K$, in the case $i=1$, without
losing generality.  The case $K=1$ is obvious, since in this case
$M=(1)$, $Q^{(1)}=(\sigma^2)$ is indeed the variance of the offspring
distribution and $\Pi^{(1)}$ is the identity.  According to Lemma
\ref{toutvabien}, under $P^\bx$, it is licit to do the $K$-to
$K-1$-type operation $\wt{\Pi}$, without changing the hypothesis that
the GW processes under consideration are nondegenerate, irreducible
and critical. This immediately gives the result on the mean of the
offspring distribution of $\Pi^{(1)}(F)$ by induction.  The only
statement that remains to be proved is the formula for the variance of
its offspring distribution.  Using again (\ref{phij}) and
(\ref{phiK}), straightforward (but tedious) computations show that,
letting
$$\wt{Q}^{(i)}_{jk}=\frac{\dd^2 \wt{\varphi}^{(i)}}{\dd s_j\dd s_k}({\bf 1})\, 
,\qquad 1\leq i,j,k\leq K$$ be the quadratic forms associated with the
offspring distributions of $\wt{\Pi}(F)$ under $P^\bx$,
\begin{eqnarray*}
  \wt{Q}^{(i)}_{jk}&=&
  Q^{(i)}_{jk}+\frac{m_{Kk}Q^{(i)}_{jK}+m_{Kj}Q^{(i)}_{kK}
  }{1-m_{KK}}+\frac{m_{Kj}m_{Kk}}{(1-m_{KK})^2}Q^{(i)}_{KK}\\
  & & +\frac{m_{iK}}{1-m_{KK}}\left(Q^{(K)}_{jk}+
    \frac{m_{Kk}Q^{(K)}_{jK}+m_{Kj}Q^{(K)}_{kK}
    }{1-m_{KK}}+\frac{m_{Kj}m_{Kk}}{(1-m_{KK})^2}Q^{(K)}_{KK}\right).
\end{eqnarray*}
It is then easy to check that
$$\wt{\ba}\cdot\wt{\bQ}(\wt{\bb})=
\frac{1-a_K}{(1-a_Kb_K)^2}\ba\cdot\bQ(\bb).$$ Using the induction
hypothesis, we obtain
$\wt{\ba}\cdot\wt{\bQ}(\wt{\bb})=\wt{a}_1\wt{b}_1^2\var(\ov{\mu}^{(1)}_1)$,
so that
$$a_1b_1^2\frac{1-a_K}{(1-a_Kb_K)^2}\var(\ov{\mu}^{(1)}_1)=
\frac{1-a_K}{(1-a_Kb_K)^2}\ba\cdot\bQ(\bb),$$ giving the result.

(ii) We again prove this in the case $i=1$. When $K=1$ there is
nothing to prove. If $K=2$, one checks that the number of type-$2$
vertices trapped between two $1$-type generations of $\Pi^{(1)}(T)$
under $P^{(1)}$ has mean $m_{12}/(1-m_{22})$ and finite variance
(resp.\ some exponential moment if $\mmu$ has some), by
differentiating (\ref{eq:13}) and (\ref{eq:14}) once. Then, one
obtains by direct computations that this is $a_2/a_1=(1-a_1)/a_1$.

So suppose $K\geq 3$.  The idea is to apply the projection operation
$\wt{\Pi}$, $K-2$ times, removing types $K,K-1,\ldots,3$ one after the
other. When this is performed, a two-type tree $\wt{\Pi}^{\circ k}(F)$
is obtained, and the number of type $2$ vertices that have only the
root as type $1$ ancestor is precisely the number of type $2$
individuals that are trapped between two generations of
$\Pi^{(1)}(F)$. By a direct inductive argument using Lemma
\ref{toutvabien} and the discussion after its proof, the mean matrix
of this contracted two-type tree $\wt{\Pi}^{\circ k}(F)$ has
$(a_1,a_2)$ as left $1$-eigenvector. In view of the $K=2$ case above,
the mean number of deleted type $2$ vertices in a generation is thus
$a_2/a_1$. By symmetry, the average number of type $i\in\{2,\ldots,
K\}$ vertices deleted in a generation of $\Pi^{(1)}(F)$ is
$a_i/a_1$. The average total number of deleted vertices is thus
$(a_2+\ldots+a_K)/a_1=(1-a_1)/a_1$, as claimed.

Finally, (iii) is obtained by applying point (iii) in Lemma
\ref{toutvabien} in a similar induction argument.  \cq

\subsection{Two exponential bounds}\label{sec:two-lemma}

Let $\zzeta$ be an ordered offspring distribution with $p_*\zzeta=\mmu$
satisfying {\rm (H)}. The following lemma allows to control the height
and number of components in a GW forest.

\begin{lmm}\label{sec:two-lemma-1}
  There exist two constants $0<C,C'<\infty$ depending only on $\mmu$,
  such that for every $n\in\N$, $\bx\in[K]^\N$ and $\eta>0$,
$$P^{\bx}\left(\max_{0\leq k\leq n}|u(n)|\geq n^{1/2+2\eta}\right)\leq C
\exp(-C'\,n^{\eta}),$$ and
$$P^{\bx}\left(\Upsilon^F_n\geq n^{1/2+2\eta}\right)\leq C \exp(-C' n^{\eta}
),$$
\end{lmm}

\proof In the case $K=1$, this is a straightforward consequence of
\cite[Lemma 13]{jfmgm05}. In the general case, it suffices to note
that under $P^\bx$, independently of $\bx$,
$$\max_{0\leq k\leq n}|u_F(k)|\leq \sum_{i\in[K]}\max_{0\leq k\leq
  n}|u_{\Pi^{(i)}(F)}(k)|\, ,$$
and 
$$\Upsilon^F_n\leq \sum_{i\in[K]}\Upsilon^{\Pi^{(i)}(F)}_n\, .$$
We may then apply the $K=1$ case to each of the forests
$\Pi^{(i)}(F),i\in[K]$, which under $P^\bx$ are critical
non-degenerate monotype GW forests by Lemma \ref{monotype}. \cq

\subsection{Convergence of types}\label{sec:convtypes}

A natural way to proceed to prove Theorem \ref{T1} is now to use the
known results of the monotype forest $\Pi^{(1)}(F)$, and to try and
pull them back to the projected multitype tree.  To do this, we must
take care of two kinds of loss of information: the number of vertices
with type $\neq 1$ of $F$ that stand between two consecutive $1$-type
vertices seen in $\Pi^{(1)}(F)$ (`time' information), and the number
of vertices of $F$ that actually stand between a type-$1$ vertex of
$\Pi^{(1)}(F)$ and one of its sons (`height' information).

For $u\in\ff\in\FF^{(K)}$, let $\anc_\ff^u(i):= \anc_\ff^u(e(v)=i)$ be
the number of ancestors $u$ such that $e(v)=i$.

\begin{prp}\label{HEIGHT}
  Under {\rm (H)}, for every $\gamma>0$ and $\bx\in[K]^\N$,
  \begin{equation}
    \max_{i\in[K]}  P^{\bx}\left(\max_{0\leq k\leq n}\left|H^F_k-
        \frac{\anc_F^{u(k)}(i)}{a_ib_i}\right|>
      n^{1/4+\gamma}\right)= \oe(n)\, .
\end{equation}
\end{prp}

For this, we need the following moderate deviations estimate for
Markov chains.

\begin{lmm}\label{sec:convergence-types}
  On some probability space $(\Omega,\a,P)$, let $(X_n,n\geq 0)$ be an
  irreducible Markov chain taking values in a finite set $S$. Let
  $\pi$ be its stationary distribution, and
  $\xi_n=n^{-1}\sum_{k=0}^{n-1}\delta_{X_k}$ be its empirical
  distribution at time $n$. Then, for any $f: S\to \R$ and $\gamma>0$,
  there exists $N(f,\gamma)>0$ for every $n\geq N(f,\gamma)$,
$$\max_{0\leq k\leq n}P(k|\xi_k(f)-\pi(f)|\geq n^{1/2+\gamma})\leq 
\exp(-n^{\gamma})\, .$$
\end{lmm}

\proof If the Markov chain is also aperiodic, then according to Wu
\cite[Theorem 2.1 (a)]{wu95}, for every $f$ there exists a constant
$C$ depending on $f$ such that for every large enough $k$,
$$
P(k|\xi_k(f)-\pi(f)|\geq k^{1/2+\gamma})\leq \exp(-Ck^{2\gamma})\, .$$ If
the chain has period $d>1$, then the same result is easily obtained by
partitioning the state space into the $d$ periodic classes with equal
$\pi$-masses and considering the $d$ shifted chains $(X_{nd+r},n\geq
0)$ for $0\leq r\leq d-1$, which under $P(\cdot|X_0=x)$ are aperiodic
for every $x\in S$.

Next, notice that for large enough $n$, $\max_{0\leq k\leq
  n^{1/2+\gamma/2}}P(k|\xi_k(f)-\pi(f)|\geq n^{1/2+\gamma})=0$ since
$|\xi(f)-\pi(f)|\leq 2\|f\|_{\infty}$. Thus, we have, for large enough
$n$,
\begin{eqnarray*}
\max_{0\leq k\leq n}P(k|\xi_k(f)-\pi(f)|\geq n^{1/2+\gamma})&\leq& 
\max_{n^{1/2+\gamma/2}\leq k\leq n}P(k|\xi_k(f)-\pi(f)|\geq
k^{1/2+\gamma})\\ 
&\leq& \exp(-Cn^{\gamma+\gamma^2})\, ,
\end{eqnarray*} 
entailing the result. \cq

\bigskip

\noindent{\bf Proof of Proposition \ref{HEIGHT}. } Suppose $i=1$ with
no loss of generality.  By Lemma \ref{sec:two-lemma-1}, the
probability that either $\Upsilon^F_{n}>n^{1/2+\gamma}$ or
$\max_{0\leq k\leq n}|u(k)|>n^{1/2+\gamma}$ is an $\oe(n)$, so we can
restrict ourselves to the complementary event. Thus, it suffices to
bound the quantity
\begin{eqnarray*}
  P^\bx\left(\ind_{\left\{\Upsilon^F_n\leq n^{1/2+\gamma},
        \max_{0\leq k\leq
          n}|u(k)|\leq n^{1/2+\gamma}\right\}}\max_{0\leq k\leq n}
    |a_1b_1|u(k)|-\anc_F^{u(k)}(1)|>
    n^{1/4+\gamma}\right)\\
  \leq  P^\bx_{\lfloor n^{1/2+\gamma}\rfloor}\left(\max_{u\in F,|u|\leq
      n^{1/2+\gamma}}|a_1b_1|u|-\anc_F^u(1)|>n^{1/4+\gamma}
  \right)
  \end{eqnarray*} 
  for large $n$. By bounding the max by a sum over the same set, and
  then making use of the ancestral decomposition (Lemma
  \ref{ancestral}), this is less than
\begin{eqnarray*}
  \lefteqn{E^\bx_{\lfloor n^{1/2+\gamma}\rfloor}
    \left[\sum_{u\in F}\ind_{\{|u|\leq n^{1/2+
          \gamma}\}}\ind_{\left\{|a_1b_1|u|-\anc_F^u(1)|>
          n^{1/4+\gamma}\right\}}\right]}\\
  &\leq&  C\sum_{j=1}^{\lfloor n^{1/2+\gamma}\rfloor} 
  \sum_{h=0}^{\lfloor n^{1/2+\gamma}\rfloor}
  \wh{P}^{(x_j),h}\left(h|h^{-1}\anc_T^V(1)-a_1b_1|>
    n^{1/4+\gamma}  \right)\\
  &\leq &Cn^{1+2\gamma}\max_{i\in[K]}\max_{0\leq h\leq
    n^{1/2+\gamma}}\wh{P}^{(i),h}\left(h|h^{-1}\anc_T^V(1)-a_1b_1|>
    n^{1/4+\gamma}  \right)\, ,
\end{eqnarray*}
where $C=\max_i b_i/\min_i b_i>0$.  Recall that under
$\wh{P}^{(i),h}$, the sequence $(e(V_0),\ldots,e(V_h)=e(V))$ is a
Markov chain in $[K]$ started at $i$ with step transition 
$p_{j,j'}=m_{jj'}b_{j'}/b_j$, and which admits
$\ba\bb:=(a_1b_1,\ldots,a_Kb_K)$ as invariant probability. Notice that
$h^{-1}\anc_T^V(1)$ is the empirical measure of $\{1\}$ for this
Markov chain. The result is now a straightforward consequence of Lemma
\ref{sec:convergence-types}. \cq

\bigskip

Next, recall the notation $\Lambda^\ff_i(k)$ and let also $G^\ff_i(k)=
\card\{u\prec u^{(i)}(k)\}$, where $u^{(i)}(0),u^{(i)}(1),\ldots$ is
the list of type-$i$ vertices of $\ff$, arranged in depth-first
order. A similar notation holds for trees instead of forests, and we
adopt the convention $G^\bt_i(\#\bt^{(i)})=\#\bt$.

\begin{prp}\label{J}
{\rm (i)}  For any $\bx\in[K]^\N$, under $P^{\bx}$, as $n\to\infty$,
$(\Lambda^F_i(\lfloor ns\rfloor)/n,s\geq 0)$ converges in probability to
$s\mapsto a_i s$, for the topology of uniform convergence over compact
sets.

{\rm (ii)} Moreover, if $\mmu$ admits small exponential moments, it
holds that for every $\gamma>0$ and $\bx\in[K]^\N$,
\begin{equation}\label{mdj}
  P^\bx\left(
    |G^F_i(n)-a_i^{-1} n|>n^{1/2+\gamma}
  \right)=\oe(n)\, .
\end{equation}
\end{prp}

\proof With the notations of Section \ref{sec:proj-monotype-trees},
for $\ff\in\FF^{(K)}$, let let $N(k):=N^{(i)}(u(k))$ be the number of
descendents $v$ of $u^{(i)}(k)$ such that the types of vertices in
$\{w:u^{(i)}(k)\vdash w\vdash v, w\neq u^{(i)}(k)\}$ are all $\neq i$,
and $N'(k)$ is the similar quantity, but counting only the vertices
$w$ that come before $u^{(i)}(n)$ in depth-first order. Then
\begin{equation}  
G^\ff_i(n)=
  \sum_{k=0}^{n-1}(1+N(k)) + \underbrace{\sum_{k=0}^{n-1}(N'(k)-N(k))
    \ind_{\{u^{(i)}(k)\vdash u^{(i)}(n)\}} }_{ :=R_1(n)}
  +\underbrace{\sum_{k=1}^{\Upsilon^\ff_n}N^{(i)}_{\rm
      floor}(k)}_{:=R_2(n)}
  \, ,\label{eq:15}
\end{equation}
Now, we estimate the
probability that $R_1(n)+R_2(n)$ is large, and by Lemma
\ref{sec:two-lemma-1}, for any fixed $0<\delta<1/2$, we may restrict
ourselves to the event that the number of ancestors of type $1$ of
$u^{(i)}(n)$ is $\leq n^{1/2+\delta}$ and that the tree containing
$u^{(i)}(n)$ has rank $\leq n^{1/2+\delta}$, up to losing an
$\oe(n)$ term.  Then
under this event, the probability of $\{R_2(n)>n^{1-\delta}\}$ is less
than
$$P^\bx\left(\sum_{k=1}^{n^{1/2+\delta}}N^{(i)}_{\rm floor}(k)>
  n^{1-\delta}\right)\, .$$ Now under $P^\bx$, the $N^{(i)}_{\rm
  floor}(k)$'s are independent, with respective laws that of
$\card\{v:i\notin\{e(w):w\vdash v\}\}$ under $P^{(x_k)}$. By (ii) in
Proposition \ref{monotype} and Chebychev's inequality, this
goes to $0$ as $n\to\infty$.  

On the other hand, notice that $|R_1(n)|\leq \sum_{k=0}^{n-1} N(k)
\ind_{\{u^{(i)}(k)\vdash u^{(i)}(n)\}}$. Let
$r_n=\sup\{k:P^\bx(N(0)>k)> n^{-1}\}$.  Since $N(0)$ has finite
variance under $P^\bx$ by (ii) in Proposition \ref{monotype}, it holds
that $P^\bx(N(0)>t)=o(t^{-2})$ as $t\to\infty$, so that
$r_n=o(n^{1/2})$, as otherwise $r_{\phi(n)}\geq c\phi(n)^{1/2}$ for
some extraction $\phi$, so that
$$\phi(n)^{-1}<P^{\bx}(N(0)>r_{\phi(n)})\leq
P^\bx(N(0)>c\phi(n)^{1/2})=o(\phi(n)^{-1})\, ,$$ a
contradiction. Therefore, for any $r'_n=o(n^{1/2})$ such that
$r_n=o(r'_n)$, we obtain that 
$$P^{\bx}\left(\max_{0\leq k\leq
    n}N(k)>r'_n\right)\leq
1-(1-P(N(0)>r'_n))^n\build\longrightarrow_{n\to\infty}^{}0\, .$$ On the
other hand, we know from (i) in Proposition \ref{monotype} and
Proposition \ref{thduqleg} that
$$n^{-1/2}\max_{0\leq k\leq n}\anc_F^{u(k)}(i)\leq
n^{-1/2}\max_{0\leq k\leq n}H^{\Pi^{(i)}(F)}_k\, ,$$ 
which converges in
distribution as $n\to\infty$ to the supremum of a properly scaled
Brownian excursion. Consequently, noticing that $R_1(n)$ is a sum
involving $\anc_F^{u^{(i)}(n)}(i)$ terms, we obtain that for every
$\eps>0$, there exists $C>0$ such that
$$P^\bx\left(n^{-1}R_1(n)\leq Cn^{-1}n^{1/2}r'_n\right)>1-\eps$$
for every $n$ large, whence $R_1(n)=o(n)$ in probability.

These estimates, when combined with (\ref{eq:15}) and the law of large
numbers, entail that 
$G^F_i(n)/n$ converges in probability to the mean of $1+N(0)$
under $P^{(i)}$, which by (ii) in Proposition \ref{monotype} is
$1/a_i$.  Therefore, $G^F_i(\lfloor ns\rfloor)/n\to a_i^{-1}s$ in probability
for every rational $s$ and we claim that the convergence holds for the
uniform topology over compact subsets of $\R$. To see this, one can
use Skorokhod's representation theorem and assume that the convergence
of $G^F_i(\lfloor ns\rfloor)$ is almost-sure for every rational $s$, and then
apply a standard monotonicity, continuity and compactness argument.
It is then elementary to conclude that the right-continuous inverse
function $(\Lambda^F_i(\lfloor ns\rfloor)/n,s\geq 0)$ converges in probability
to $s\mapsto a_i s$ for the uniform topology over compact sets.

Part (ii) of the statement is obtained along closely related lines, by
first noting that this time, for $0<\delta<\gamma/2$, and using
similar notations as above,
\begin{eqnarray*}
  P^\bx\left(
R_2(n)>n^{1/2+2\delta}\right) & \leq & 
P^\bx\left(\sum_{k=1}^{n^{1/2+\delta}}N^{(i)}_{\rm floor}(k)>
  n^{1/2+2\delta}\right)+
\oe(n)\\
&\leq &
\exp(An^{1/2+\delta}-\eps n^{1/2+2\delta})+\oe(n)=\oe(n)
\end{eqnarray*}
for some $A,\eps>0$. Then, we have $P^\bx\left(\max_{0\leq k\leq n}
  N(k)>n^\delta\right)=\oe(n)$, so that
$$P^\bx\left(
  R_1(n)>n^{1/2+2\delta}\right)\leq P^\bx\left(\max_{0\leq k\leq
    n}|u(k)|\geq n^{1/2+\delta}\right)+ P^\bx\left(\max_{0\leq k\leq
    n}N(k)>n^{\delta}\right)=\oe(n)\, .$$ Finally, the estimate
$$P^\bx\left(
  \left|\sum_{k=0}^{n-1}(1+N(k))-a_i^{-1}n\right|\geq
  n^{1/2+\gamma}\right)=\oe(n)$$ is a standard moderate deviations
estimate for random variables admitting small exponential moments, see
\cite[Theorem 2.6]{petrov95}. This is enough to conclude.  \cq

\bigskip

Notice that the previous statement immediately implies point (iii) in
the statement of Theorem \ref{T1}.

\subsection{Proof of (i) and (ii) in Theorem
  \ref{T1}}\label{sec:proof-i-ii}

For any $s\geq 0$, we have
$$\left|H^F_{\lfloor ns\rfloor}-
  \frac{H^{\Pi^{(i)}(F)}_{\Lambda^F_i(\lfloor ns\rfloor)-1}}{
    a_ib_i}\right|\leq \left|u(\lfloor
  ns\rfloor)-\frac{\anc_F^{u(\lfloor ns\rfloor)}(i)}{a_ib_i}\right|
+\frac{\left|H^{\Pi^{(i)}(F)}_{\Lambda^F_i(\lfloor
      ns\rfloor)-1}-\anc_F^{u(\lfloor ns\rfloor)}(i)
  \right|}{a_ib_i}\, .$$ By Proposition \ref{HEIGHT}, we obtain that,
for every $A>0$,
$$\sup_{0\leq s\leq A}n^{-1/2}\left|u(\lfloor ns\rfloor)-
  \frac{\anc_F^{u(\lfloor ns\rfloor)}(i)}{a_ib_i}\right|\to 0 $$ in
probability as $n\to\infty$. On the other hand, we claim that
\begin{equation}\label{eq:8}
  \left|H^{\Pi^{(i)}(F)}_{\Lambda^F_i(\lfloor ns\rfloor)-1}-
    \anc_F^{u(\lfloor ns\rfloor)}(i)
  \right|\leq \left|H^{\Pi^{(i)}(F)}_{\Lambda^F_i(\lfloor ns\rfloor)-1}-
    H^{\Pi^{(i)}(F)}_{\Lambda^F_i(\lfloor ns\rfloor)}\right|+1\, .
\end{equation} 
Indeed, if $u^{(i)}(\lfloor ns\rfloor)$ is an ancestor of $u(\lfloor ns\rfloor)$, then the
left-hand side is zero and there is nothing to prove. Else, the
left-hand side equals the number of ancestors of type $i$ of
$u^{(i)}(\lfloor ns\rfloor)$ which are not ancestors of $u(\lfloor ns\rfloor)$, and so
$H^{\Pi^{(i)}(F)}_{\Lambda^F_i(\lfloor ns\rfloor)-1}-\anc_F^{u(\lfloor ns\rfloor)}(i)\geq 0$. On
the other hand, the strict ancestors of $u^{(i)}(\lfloor ns\rfloor+1)$ that are not
ancestors of type $i$ of $u(\lfloor ns\rfloor)$, cannot be themselves of type $i$
by definition (otherwise, such an ancestor would come after $u(\lfloor ns\rfloor)$
and before $u^{(i)}(\lfloor ns\rfloor+1)$ in depth-first order). Hence,
$H^{\Pi^{(i)}(F)}_{\Lambda^F_i(\lfloor ns\rfloor)}-\anc_F^{u(\lfloor ns\rfloor)}(i)\leq 1$, so
that 
$$0\leq H^{\Pi^{(i)}(F)}_{\Lambda^F_i(\lfloor ns\rfloor)-1}-\anc_F^{u(\lfloor ns\rfloor)}(i)
\leq H^{\Pi^{(i)}(F)}_{\Lambda^F_i(\lfloor ns\rfloor)-1}-
  H^{\Pi^{(i)}(F)}_{\Lambda^F_i(\lfloor ns\rfloor)}+1\, ,$$ and the
claimed inequality follows. 

Under $P^{\bx}$, the forest $\Pi^{(i)}(F)$ is a single-type GW forest
whose offspring distribution has finite variance by Proposition
\ref{monotype}, so that by Proposition \ref{thduqleg}, 
$$n^{-1/2}\max_{0\leq k\leq n}|H^{\Pi^{(i)}(F)}_{k-1}-
H^{\Pi^{(i)}(F)}_k|\build\longrightarrow_{n\to\infty}^{P^{\bx}}0\, ,$$
and it follows that under $P^{\bx}$,
\begin{equation}\label{eq:23}
  \left(n^{-1/2}\left(H^F_{\lfloor ns\rfloor}-(a_ib_i)^{-1}H^{\Pi^{(i)}(F)}_{
        \Lambda^F_i(\lfloor ns\rfloor)}\right),s\geq 0\right)
  \build\longrightarrow_{n\to\infty}^{P^\bx}0 
\end{equation}
for the topology of uniform convergence over compact
sets.

Using Propositions \ref{J}, \ref{monotype} and \ref{thduqleg}, and
composing $s\mapsto H^{\Pi^{(i)}(F)}_{\lfloor ns\rfloor}$ with $s\mapsto
\Lambda^F_i(\lfloor ns\rfloor)/n$, we now obtain that
$(n^{-1/2}H^{\Pi^{(i)}(F)}_{\Lambda^F_i(\lfloor ns\rfloor)},s\geq 0)$ converges in
distribution to $(2\ov{\sigma}_i^{-1}|B_{a_is}|,s\geq 0)$ where
$\ov{\sigma}_i^2=\var(\ov{\mu}^{(i)}_i)$, which is also equal in law
to $(2a_i^{1/2}\ov{\sigma}_i^{-1}|B_s|,s\geq 0)$.  One way of seeing
this is to use Skorokhod's representation theorem to exhibit a
probability space where the convergences of Propositions
\ref{thduqleg} and \ref{J} hold a.s.\ rather than in distribution.
Point (i) of the theorem is now proved by using (\ref{eq:23}).

Let us prove (ii).  By definition, since all the roots are of type
$i$, $u(\lfloor ns\rfloor)$ and the last node with type $i$ before $u(\lfloor ns\rfloor)$ in
depth-first order belong to the same tree. Therefore, the label of the
tree of $F$ containing $u(\lfloor ns\rfloor)$ is always the same as the label of
the tree of $\Pi^{(i)}(F)$ containing the $\Lambda^F_i(\lfloor ns\rfloor)$-th
node. This implies that
$\Upsilon^{\Pi^{(i)}(F)}_{\Lambda^F_i(\lfloor ns\rfloor)}=\Upsilon^{F}_{\lfloor ns\rfloor}$. Now,
the result is a plain consequence of Proposition \ref{thduqleg} and of
similar arguments as above. \cq

\bigskip

\noindent{\bf Proof of Corollary \ref{sec:main-result}. }
From (ii) in Theorem \ref{T1}, we obtain that
$(n^{-1}H^{F}_{n^2s},0\leq s\leq \tau_n)$ converges in
distribution to $2\sigma^{-1}(|B_s|,0\leq s\leq
\mathfrak{T}_{b_i\sigma^{-1}})$, where $\tau_n$ is the first
hitting time of $n$ by $\Upsilon^{F}$ and $\mathfrak{T}_{x}$ is the
first hitting time of $x$ by $L^0$.

Now, 
$$
  P^{{\bf i}}\left(\forall 1\leq k\leq n, \HT(F_k)<n\right) 
  \build\longrightarrow_{n\to\infty}^{} P(2\sigma^{-1}|B_s|
  \leq 1,0\leq s\leq 
  \mathfrak{T}_{b_i\sigma^{-1}})\, ,$$
which can be rewritten as
$$
  \left(1-P^{(i)}(\HT(T)\geq n)\right)^n \build\longrightarrow_{
    n\to\infty}^{} 
  \exp\left(-\frac{b_i}{\sigma}N\left(\frac{2}{\sigma}\sup{\rm e} \geq 1\right)
  \right),
$$
where $N(\d{\rm e})$ is the Ito excursion measure of the standard
Brownian motion (see e.g.\ \cite[Chapter XII]{revyor} for definitions
and the results recalled below), and where we have used the Ito
decomposition of a Brownian motion into a Poisson process of
excursions in the local time scale. Taking logarithms and using
$N(\sup {\rm e}\geq x)=1/x$, gives the result. \cq

Let us also mention that similar arguments, following the same lines
as in \cite[Proposition 2.5.2]{duqleg02}, actually show the more
general result:

\begin{crl}\label{sec:proof-theorem-reft1}
  For every $a>0$ the probability measures $P^{(i)}(n^{-1}H^{T}\in
  \cdot\, |\HT(T)\geq an)$ converge in distribution as $n\to\infty$
  towards $N(2\sigma^{-1}{\rm e}\in \cdot\,| 2\sigma^{-1}\sup {\rm
    e}\geq a)$.
\end{crl}

\subsection{Conditioned results: Theorem
  \ref{sec:convergence-gw-trees}}\label{sec:conditioned-results}

Our main tool for conditioning is the following estimate for the size
of GW trees. 

\begin{lmm}\label{sec:cond-results:-theor-1}
Let $\mmu$ be a critical non-degenerate offspring distribution with
finite variance. Then for every $i,j\in[K]$, one has 
$$n^{3/2}P^{(i)}(\#T^{(j)}=n)\build\longrightarrow_{n\to\infty}^{}C_{ij}\, ,$$
where if necessary the limit is taken along a subsequence for which
the probability on the left-hand side is non-zero, and for some
constant $C_{ij}>0$.  
\end{lmm}

\proof This is very similar to Lemma 14 in \cite{jfmgm05}. If $i=j$,
then using the fact that the reduced tree of Section
\ref{sec:proj-monotype-trees} is a monotype GW tree, the result is a
well-known fact.  To treat the general case we elaborate slightly on
the proof.

Let $r\in\N$ be fixed, and recall ${\bf j}=(j,j,\ldots)$ and the
notation at the very end of Section \ref{sec:galton-watson-trees}.
Let $\tilde{\mu}:=\ov{\mu}^{(j)}_j$ be the offspring distribution of
the GW tree $\Pi^{(j)}(T)$ under $P^{(j)}$, and on some probability
space $(\tilde{\Omega},\tilde{\a},\tilde{P})$, let $(W_n,\geq 0)$ be a
random walk with step distribution $\tilde{\mu}(\cdot+1)$ on
$\{-1,0,1,2,\ldots\}$. Under $P^{{\bf j}}_r$, the forest
$\Pi^{(j)}(F)$ is a monotype GW forest with offspring distribution
$\tilde{\mu}$ and $r$ tree components.  It is then well-known that
$$P^{{\bf j}}_r(\#\Pi^{(j)}(F)=n)=P^{{\bf
    j}}_r(\#F^{(j)}=n)=\frac{r}{n}\tilde{P}(W_n=-r)\, .$$ By the local
limit theorem in the lattice case \cite[Theorem XV.5.3]{feller2},
$n^{1/2}\tilde{P}(W_n=-r)\to C$ as $n\to\infty$ for some $C>0$, and
for every $r$ such that the probabilities under consideration are
$>0$. Moreover, there is a common uniform bound for all the terms as
$r$ varies along the admissible values. 

Let $p(r;i,j)$ be the probability that there are $r$ tree components
in $\Pi^{(j)}(T)$ under $P^{(i)}$. Notice that the probability
distribution $(p(r;i,j),r\in\N)$ has finite expectation (its
generating function is $\ov{\varphi}_i^{(j)}$ with the notations of
Section \ref{sec:proj-monotype-trees}), so that $\sum_r
rp(r;i,j)<\infty$. Then 
$$P^{(i)}(\#T^{(j)}=n)=\sum_{r\geq 1}p(r;i,j)P^{{\bf
    j}}_r(\#\Pi^{(j)}(F)=n)\, ,$$ and an application of the previous
paragraph and dominated convergence (using the fact that $n^{3/2}P^{{\bf
    j}}_r(\#\Pi^{(j)}(F)=n)$ is uniformly bounded) gives that
\begin{equation}\label{eq:24}
n^{3/2}P^{(i)}(\#T^{(j)}=n)\to C\sum_{r\geq 1}rp(r;i,j)\, ,
\end{equation}
which is the wanted result. \cq

\begin{lmm}\label{sec:cond-results:-theor-4}
  The respective laws of the number of tree components of
  $\Pi^{(j)}(T)$ under the probability distributions
  $P^{(i)}(\cdot|\#T^{(j)}=n)$ converge weakly as $n\to\infty$.
\end{lmm}

\proof We use the notations of the previous proof, as well as the
expression (\ref{eq:24}) of the constant $C_{ij}$. Observe that the
$P^{(i)}$-probability that $\Pi^{(j)}(T)$ has $r$ components given it
has $n$ individuals is
$$P^{(i)}(c_{\Pi^{(j)}(T)}(\varnothing)=r|
\#T^{(j)}=n)=\frac{p(r;i,j)P^{{\bf
      j}}_r(\#F^{(j)}=n)}{P^{(i)}(\#T^{(j)}=n)}\build
\longrightarrow_{n\to\infty}^{}
\frac{rp(r;i,j)}{\sum_{r'}r'p(r';i,j)}\, ,$$ and this does define a
probability distribution. \cq

The following modification of Proposition \ref{thduqleg} for forests
with a fixed number of trees also holds:

\begin{lmm}\label{sec:cond-results:-theor-2}
  In the case $K=1$, assume {\rm (H)} and take $r\in\N$. Let
  $P_r:=P^{{\bf 1}}_r$ be the law of a
  monotype GW forest with $r$ tree components and offspring
  distribution $\mu:=\mu^{(1)}$. Then the process $(n^{-1/2}H^F_{\lfloor n
    s\rfloor},0\leq s\leq 1)$ under $P_r(\cdot|\#F=n)$ converges in
  distribution to $2\sigma^{-1}B^{\rm ex}$ as $n\to\infty$.
\end{lmm}

\proof The law of the total size of a GW tree is in the domain of
attraction of a positive stable law with index $1/2$, as follows from
the previous lemma in the case $K=1$, so that when taking $r$
independent copies of GW trees with offspring distribution $\mu$, and
conditioning their sum to be $n$, only one of the trees has a size of
order $n$, while the others $r-1$ trees have $o(n)$ size, hence have
maximal height $o(n^{1/2})$ according to Lemma
\ref{sec:two-lemma-1}. Hence the result.  \cq

\bigskip

\noindent{\bf Proof of Theorem \ref{sec:convergence-gw-trees}. }
The proof starts like the one of Theorem \ref{T1}. For $0\leq s\leq
1$, write
\begin{equation}\label{eq:16}
  \left|H^T_{\lfloor\#T s\rfloor}-\frac{H^{\Pi^{(j)}(T)}_{
        \Lambda^T_j(\lfloor \#T s\rfloor)}}{a_jb_j}
  \right|\leq \left|u(\lfloor \#T s\rfloor)-
    \frac{\anc_T^{u(\lfloor \#T s\rfloor)}(j)}{a_jb_j}\right|
  +R_n(s)\, ,
\end{equation} 
where $|R_n(s)|\leq (a_jb_j)^{-1}(2\max_{0\leq k\leq
  n}|H^{\Pi^{(j)}(T)}_{k-1}-H^{\Pi^{(j)}(T)}_k|+1)$.

We start by showing the convergence of the processes
$(n^{-1/2}H^{\Pi^{(j)}(T)}_{\Lambda^T_j(\lfloor \#T s\rfloor)}
/a_jb_j,0\leq s\leq 1)$ under $P^{(i)}(\cdot|\#T^{(j)}=n)$. Since
$\Pi^{(j)}(T)$ under $P^{(i)}$ is a GW forest, and by first
conditioning on $c_{\Pi^{(j)}(T)}(\varnothing)=r$, we obtain using
Lemmas \ref{sec:cond-results:-theor-4} and
\ref{sec:cond-results:-theor-2} that
$(n^{-1/2}H^{\Pi^{(j)}(T)}_{\lfloor ns\rfloor},0\leq s\leq 1)$ under
$P^{(i)}(\cdot|\#T^{(j)}=n)$ converges in distribution to
$2\ov{\sigma}_j^{-1}B^{\rm ex}$, where
$\ov{\sigma}_j^2=\var(\ov{\mu}^{(j)}_j)$, with definitions from
Proposition \ref{monotype}.

We now show that $(n^{-1}\Lambda^T_j(\lfloor \#T s\rfloor),0\leq s\leq 1)$
converges in probability to the identity on $[0,1]$, under the
conditioned measures. Using Lemma \ref{sec:cond-results:-theor-1} and
(ii) in Proposition \ref{J}, one obtains that for some $C>0$, for
every $s\in[0,1]$,
\begin{eqnarray}
 \lefteqn{
P^{{\bf i}}\left(\left|G^F_j(\lfloor ns\rfloor)-a_j^{-1}sn
      \right|\geq n^{1/2+\gamma}\big| \#F_1^{(j)}=n\right)}
\nonumber\\
  &
\leq
& 
C 
  n^{3/2}P^{{\bf i}}\left(\left|G^F_j(\lfloor ns\rfloor)-a_j^{-1}sn
    \right|\geq n^{1/2+\gamma}\right)
=\oe(n)\label{eq:22}
\end{eqnarray}
Since $F_1$ under $P^{\bf i}$ has same distribution as $T$ under
$P^{(i)}$, 
we thus obtain that for $0\leq s\leq 1$, 
\begin{equation}\label{eq:17}
P^{(i)}(|G^T_j(\lfloor ns\rfloor)-a_j^{-1}sn|\geq
n^{1/2+\gamma}|\#T^{(j)}=n)=\oe(n)\, ,
\end{equation} 
which shows that $n^{-1}G^T_j(\lfloor ns\rfloor)$ under
$P^{(i)}(\cdot|\#T^{(j)}=n)$ converges to $a_j^{-1}s$ for every
$s\in[0,1]$ in probability, and thus, by the same reasoning as in
Section \ref{sec:proof-i-ii}, we obtain that $(n^{-1}G^T_j(\lfloor ns\rfloor),0\leq
s\leq 1 )$ converges in probability to $(a_j^{-1}s,0\leq s\leq 1)$ for
the uniform topology In particular, for $s=1$ we obtain that
$n^{-1}\#T$ under $P^{(i)}(\cdot|\#T^{(j)}=n)$ converges to $a_j^{-1}$
in probability, where we recall that we adopted the convention
$G^T_j(n)=\#T$.

Now, $(n^{-1}\Lambda^T_j(\lfloor \#T s\rfloor),0\leq s\leq 1)$ is the
right-continuous inverse of $(\#T^{-1}G^T_j(\lfloor ns\rfloor),0\leq s\leq 1)$,
and as such, it converges to the identity of $[0,1]$ in probability
for the uniform topology. 

It remains to show that the two terms on the right-hand side of
(\ref{eq:16}) are $o(n^{1/2})$ in probability, uniformly in
$s\in[0,1]$.  First, notice that letting
$\Upsilon^j:=c_{\Pi^{(j)}(T)}(\varnothing)$ be the number of tree
components of $\Pi^{(j)}(T)$, then the law of $\Pi^{(j)}(T)$ under
$P^{(i)}(\cdot| \Upsilon^j=r)$ is the same as that of $\Pi^{(j)}(T)$
under $P^{{\bf j}}_r$, i.e.\ is that of a monotype GW forest with $r$
tree components.  Using Lemma \ref{sec:cond-results:-theor-2}, one
concludes that $P^{(i)}(\sup_{0\leq s\leq 1}n^{-1/2}|R_n(s)|\geq
\eps|\#T^{(j)}=n,\Upsilon^j=r)$ converges to $0$ for any
$\eps>0$. Using Lemma \ref{sec:cond-results:-theor-1}, we know that
the laws of $\Upsilon^j$ under $P^{(i)}(\cdot|\#T^{(j)}=n)$ are tight
as $n$ varies, so that we get $P^{(i)}(\sup_{0\leq s\leq
  1}n^{-1/2}|R_n(s)|\geq \eps|\#T^{(j)}=n)\to 0$ as well.

Finally, by applying (\ref{eq:17}) for $s=1$, we obtain that
$P^{(i)}(\#T>An|\#T^{(j)}=n)=\oe(n)$ for any $A>a_j^{-1}$. Combining that
with Proposition \ref{HEIGHT}, gives for such $A$ and some $C>0$:
\begin{eqnarray*}
  \lefteqn{P^{(i)}\left(\max_{0\leq k\leq
        \#T}\left|u(k)-\frac{\anc_T^{u(k)}(j)}{a_jb_j}
      \right|\geq n^{3/8}\big|\#T^{(j)}=n\right)}\\
  &\leq & Cn^{3/2}
  P^{{\bf i}}\left(\max_{0\leq k\leq An}\left|u(k)-
      \frac{\anc_F^{u(k)}(j)}{a_jb_j}
    \right|\geq n^{3/8}\right)+\oe(n)=\oe(n)\, ,
\end{eqnarray*}
hence the result. \cq

\bigskip

\rem In the companion paper \cite{mierinv}, a similar statement as
Theorem \ref{sec:conv-mult-snak} was needed, with the law
$P^{(i)}(\cdot|\#T^{(j)}=n)$ replaced by
$P^{(i,i)}(\cdot|\#F^{(j)}=n)$, i.e.\ by a forest with two trees
conditioned by their total number of vertices of type $j$. The proof
of such a statement should be clear from the previous methods: by
applying the transformation $\Pi^{(j)}$ to this forest, one obtains a
monotype Galton-Watson forest with a random number of roots that is
tight as $n$ varies, and conditioned to have $n$ vertices. By
conditioning on its number of roots and applying Lemma
\ref{sec:cond-results:-theor-2}, this implies that none but one of
these trees have more that $o(n)$ vertices. Hence in the initial
forest with two components, only one of the components has more than
$o(n)$ vertices, and the result is now a consequence of Theorem
\ref{sec:conv-mult-snak}, which will be proved in the next section.

\section{Proof of Theorem \ref{T2}}

The key technical results needed to prove Theorem \ref{T2} are, like
in \cite{jfmgm05}, a control on the frequencies of branching events in
GW trees, which will allow to prove the convergence of
finite-dimensional marginals of the snake, and a bound on a H\"older
norm-like quantity for the rescaled height process, which will imply
the tightness.

\subsection{Exponential control of branching
  events}\label{sec:expon-contr-branch}

For $u\in\ff\in\FF^{(K)}$, let
$\anc_\ff^u(i,\bw,l,h):=\anc_\ff^u(e(v)=i,\bw_\ff(v)=\bw,u\in
vl\ff_{vl},|v|\geq |u|-h)$, i.e.\ the number of ancestors $v$ of $u$
with type $i$, with children's types $\bw$, such that $u$ is the
descendent of the $l$-th child of $v$, and that are at distance at
most $h$ from $u$. In the sequel, when dealing with quantities of the
form $\anc_\ff^u({\sf P}_v)$, and if the last argument is $h$, we will
understand that we consider only those ancestors $v$ of $u$ such that
$|v|\geq |u|-h$.

\begin{lmm}\label{sec:proof-theorem-reft2-1}
  Assume $\mmu$ satisfies {\rm (H)} and has small exponential
  moments. Then for every $\bx\in[K]^\N$ and $\gamma>0$,
$$P^\bx\left(\max_{0\leq k\leq n}\max_{n^\gamma\leq h\leq |u(k)|,
\bw\in\W_K, 1\leq l\leq |\bw|}
  \frac{|\anc_F^{u(k)}(i,\bw,l,h)-h
    a_ib_{w_l}\zeta^{(i)}(\bw)|}{h^{1/2+\gamma}}\geq
  1\right)= \oe(n)\, .$$
\end{lmm}

We first state an intermediate lemma. Recall the construction of the
size-biased infinite tree $\wh{\mathfrak{t}}$ of Sect.\
\ref{sec:ancestral}, and the spinal path
$\varnothing=V_0,V_1,\ldots$. We assume $\wh{\mathfrak{t}}$ to be
constructed on some probability space $(\Omega,\a,P)$. If
$u\in\wh{\mathfrak{t}}$, we let $\bw_{\wh{\mathfrak{t}}}(u)\in \W_K$ be
the ordered sequence of its children's types, as for finite trees.

\begin{lmm}\label{sec:expon-contr-branch-1}

  \noindent{\rm (i)} The sequence $((e(V_n),e(V_{n+1})),n\geq 0)$ is a
  Markov chain with step transition $p_{(i,j),(i',j')}=\delta_{j
    i'}b_{j'}m_{jj'}/b_j $ and equilibrium measure
  $\pi_{(i,j)}=a_ib_jm_{ij}$.

\noindent{\rm (ii)}
Conditionally on $(e(V_k),k\geq 0)$, the variables
$(\bw_{\wh{\mathfrak{t}}}(V_n),V_{n+1}(n+1),n\geq 0)$ are independent,
(here $V_{n+1}(n+1)$ is the last letter of $V_{n+1}$, so that
$V_{n+1}=V_n V_{n+1}(n+1)$), with law defined by
$$P\left((\bw_{\wh{\mathfrak{t}}}(V_n),V_{n+1}(n+1))=(\bw,l)\, |\, 
  (e(V_k),k\geq
  0)\right)=\frac{\zeta^{(e(V_n))}(\bw)}{m_{e(V_n)e(V_{n+1})}}\, ,$$
for every $\bw\in\W_K$ and $1\leq l\leq |\bw|$ such that
$w_l=e(V_{n+1})$.
\end{lmm}

\proof The first statement is immediate to check since we already know
that $(e(V_n),n\geq 0)$ is Markov with step transition
$p_{jj'}=b_{j'}m_{jj'}/b_j$. 
The conditional independence of
$\bw_{\wh{\mathfrak{t}}}(V_n),V_{n+1}(n+1),n\geq 0$ is easy from the
construction of $\wh{\mathfrak{t}}$, and we have
$$P((\bw_{\wh{\mathfrak{t}}}(V_n),V_{n+1}(n+1))=(\bw,l)|e(V_k),k\geq
0)
=\frac{P((\bw_{\wh{\mathfrak{t}}}(V_n),V_{n+1}(n+1))=(\bw,l)|e(V_n))}{
  P(e(V_{n+1})=w_l|e(V_n))}\, ,$$
which amounts to the desired result. \cq

\bigskip

\noindent{\bf Proof of Lemma \ref{sec:proof-theorem-reft2-1}. }
Fix $\gamma>0$ and choose $0<\eta<\gamma^2$.  First, we claim
that 
$$P^\bx\left(\max_{0\leq k\leq n}c_F(u(k))\geq n^{\eta}\right)
=\oe(n)\, .$$ 
Indeed, under $P^\bx$, the sequences $(\bw_F(u^{(i)}(k)),k\geq 0)$,
for $i\in[K]$, are independent i.i.d.\ sequences with respective
common distribution $\zeta^{(i)}$, as follows from the Markov
branching property of \cite{MR1014449} and the fact that, when
exploring the forest in depth-first order, no information on the set
of children of the vertex explored at step $n$ can be obtained before
this step. Hence, the fact that $\mmu$ admits small exponential
moments gives that
\begin{eqnarray*}
P^\bx\left(\max_{i\in[K]}\max_{0\leq k\leq
    n} c_F(u^{(i)}(k))>n^{\eta}\right)&\leq&
Kn\max_{i\in[K]}\zeta^{(i)}(|\bw|>n^{\eta})\\
&=&
Kn\max_{i\in[K]}\mu^{(i)}(|\bz|_1>n^{\eta})\, ,
\end{eqnarray*} which by
(\ref{eq:7}) and Markov's inequality is an $\oe(n)$. Since
$$\max_{0\leq k\leq n}c_\ff(u(k))\leq \max_{i\in[K]}\max_{0\leq
  k\leq n}c_\ff(u^{(i)}(k))\, ,$$ this gives the claim.
Considering this, and using Lemma \ref{sec:two-lemma-1}, we may
restrict ourselves to showing that
\begin{equation}\label{eq:9}
  P^\bx_{\lfloor n^{1/2+\gamma}\rfloor}
  \left(\max_{u\in F,|u|\leq n^{1/2+\gamma}}\max_{(h,\bw,l)\in B_{n,|u|}}
    \frac{|\anc_F^u(i,\bw,l,h)-h
      a_ib_{w_l}\zeta^{(i)}(\bw)|}{h^{1/2+\gamma}}\geq
    1\right)= \oe(n)\, ,
\end{equation}
where $B_{n,h'}=\{(h,\bw,l)\in \Z_+\times \W_K\times \Z_+:
n^\gamma\leq h\leq h', 1\leq l\leq |\bw|\leq n^{\eta}\}$. By using
Lemma \ref{ancestral} in a similar way as in the proof of Proposition
\ref{HEIGHT}, the probability in the left-hand side of (\ref{eq:9}) is
bounded above by
\begin{equation}\label{eq:12}
C\sum_{r=1}^{\lfloor n^{1/2+\gamma}\rfloor}\sum_{h'=\lfloor n^\gamma\rfloor}^{\lfloor n^{1/2+\gamma}\rfloor}
\wh{P}^{(x_r),h'}\left(\max_{(h,\bw,l)\in
    B_{n,h'}}\frac{|\anc_T^V(i,\bw,l,h)-ha_ib_{w_l}
    \zeta^{(i)}(\bw)|}{h^{1/2+\gamma}}\geq 1\right)\, ,
\end{equation}
where $C=\max_j b_j/\min_j b_j$. Writing $Q_{i,\bw,l,h}$ for the
quotient appearing in the probability,  we have
\begin{equation}\label{eq:11}
\wh{P}^{(x_r),h'}\left(\max_{(h,\bw,l)\in B_{n,h'}}
  Q_{i,\bw,l,h} \geq 1\right)\leq n^{1/2+\gamma+\eta}
K^{n^{\eta}}\!\!\!\max_{(h,\bw,l)\in B_{n,h'}}\wh{P}^{(x_r),h'}
\left(Q_{i,\bw,l,h} \geq 1\right)\, ,
\end{equation} where
$n^{1/2+\gamma+\eta}K^{n^{\eta}}$ bounds the cardinality of $B_{n,h'}$
with $h'\leq n^{1/2+\gamma}$. Let
$\anc_\bt^u(i,j,h)=\anc_\bt^u(e(v)=i,e(vu_{|v|+1})=j,|v|\geq |u|-h)$
be the number of strict ancestors $v$ of $u\in\bt\in\TT^{(K)}$ such 
that $e(v)=i$,
$|v|\geq |u|-h$ and $e(v')=j$ whenever $v'$ is the child of $v$ with
$v\vdash v'\vdash u$, so that according to (ii) in Lemma
\ref{sec:expon-contr-branch-1}, under $\wh{P}^{(x_r),h'}$ and given
$\anc_T^V(i,j,h)=k$, the random variable $\anc_T^V(i,\bw,l,h)$ is
binomial with parameters $(k,\zeta^{(i)}(\bw)/m_{ij})$, for every
$\bw,l$ with $w_l=j$.

By Lemmas \ref{sec:expon-contr-branch-1}, (i) and
\ref{sec:convergence-types}, we have for every $n^\gamma\leq h\leq h'$,
\begin{eqnarray}
  \max_{i,j\in[K]}\wh{P}^{(x_r),h'}(|\anc_T^V(i,j,h)-ha_ib_jm_{ij}|\geq
  h^{1/2+\gamma})&\leq& \exp(-h^{\gamma})
  \nonumber\\
  &\leq& \exp(-n^{\gamma^2}) \, ,\label{eq:10}
\end{eqnarray} 
for every large enough $n$. In particular, for some $c>0$,
$\inf_{i,j\in[K]}\anc_T^V(i,j,h)\geq ch$ with
$\wh{P}^{(x_r),h'}$-probability $\leq \exp(-n^{\gamma^2})$ for every
$n^\gamma\leq h\leq h'$ and large enough $n$. Moreover,
\begin{eqnarray*}
  \lefteqn{ \wh{P}^{(x_r),h'}(Q_{i,\bw,l,h}\geq 1)}\\
  &\leq& \wh{P}^{(x_r),h'}
  \left(\left|\anc_T^V(i,\bw,l,h)-\frac{\zeta^{(i)}(\bw)
        \anc_T^V(i,w_l,h)}{m_{ij}}\right|\geq 
    \frac{h^{1/2+\gamma}}{2}\right)\\
 &  &+\wh{P}^{(x_r),h'}
  \left(|\anc_T^V(i,w_l,h)-ha_jb_{w_l}m_{ij}|\zeta^{(i)}(\bw)/m_{ij}\geq 
    \frac{h^{1/2+\gamma}}{2}\right)\, .
\end{eqnarray*}
Since $\zeta^{(i)}(\bw)\leq 1$ and using (\ref{eq:10}), the second
probability is $\leq \exp(-n^{\gamma^2})$ for large $n$, as
long as $(h,\bw,l)\in B_{n,h'}$ for $h'\geq n^\gamma$. 

Also, according to H\oe ffding's inequality for binomial
distributions, it holds that
\begin{eqnarray*}
  \wh{P}^{(x_r),h'}\left(\frac{|\anc_T^V(i,\bw,l,h)-k\zeta^{(i)}
      (\bw)/m_{ij}|}{
      h^{1/2+\gamma}} \geq 1\, \Big|\, \anc_T^V(i,j,h)=k\right)\\
  \leq 
  2\exp(-2 h^{1+2\gamma}/k)\, ,
\end{eqnarray*}
which for $h\geq n^\gamma$ and $k\leq ch$ is less than
$2\exp(-2n^{2\gamma^2}/c)$. This is enough to conclude that for every
$n^\gamma\leq h\leq h'\leq n^{1/2+\gamma}$ and large enough $n$,
$$\max_{(h,\bw,l)\in B_{n,h'}}\wh{P}^{(x_r),h'}(Q_{i,\bw,l,h}\geq 1)\leq
\exp(-n^{\gamma^2})\, .$$ Combining with (\ref{eq:11}) and then
(\ref{eq:12}), gives the result (notice that by our choice of $\eta$,
we have $\gamma^2>\eta$ so that the quantity on the right-hand
side of (\ref{eq:11}) is an $\oe(n)$).  \cq

\bigskip

We also have the following controls

\begin{lmm}\label{sec:expon-contr-branch-2}
Under the same hypotheses as Lemma \ref{sec:proof-theorem-reft2-1}, 

\noindent {\rm (i)} 
For every $\Gamma,\eta>0$ there exists $C>0$ so that for every large
$n$,
$$P^\bx\left(\max_{0\leq k\leq n}c_F(u(k))\geq C\log n\right)\leq
n^{-\Gamma}\, .$$ 

\noindent {\rm (ii)} 
For every $c>0$ there exists $\gamma\in(0,1)$ so that
$$P^\bx\left(\max_{0\leq k\leq n}\max_{n^\eta\leq h\leq |u(k)|}
  h^{-\gamma}\anc_F^{u(k)}(c_F(v)\geq c\log n,h)
\geq
1\right)=\oe(n)\, .$$

\noindent {\rm (iii)} 
For every $M>0$, there is a constant $C_M>0$ such that 
$$P^\bx\left(\max_{0\leq k\leq n}\max_{n^\eta\leq h\leq |u(k)|}
  \sup_{m\geq 1}m^Mh^{-1}\anc_F^{u(k)}(c_F(v)=m,h) \geq
  C_M\right)={\rm op}(n)\, .$$ Here $x_n={\rm op}(n)$ means that for every
$\Gamma>0$, $x_n\leq n^{-\Gamma}$ for every large $n$.
\end{lmm}

\proof Point (i) follows similar lines as the beginning of the proof
of Lemma \ref{sec:proof-theorem-reft2-1}. We bound the probability of
interest by 
$$Kn\max_{i\in[K]}\mu^{(i)}(|\bz|_1\geq C\log n)\leq
Kn^{1-C\eps}\max_{i\in[K]}\mu^{(i)}(\exp(\eps|\bz|_1))$$ for some
appropriate $\eps>0$, and choose $C$ so that $1-C\eps<-\Gamma$.

For (ii), we see that, using Lemma \ref{sec:two-lemma-1}, it suffices
to bound the probability 
$$P^\bx_{\lfloor n^{1/2+\eta}\rfloor}\left(\max_{u\in F,|u|\leq n^{1/2+\eta}}
  \max_{n^\eta\leq h\leq |u|}h^{-\gamma} \anc_F^u(c_F(v)\geq c\log
  n,h)\geq 1\right)\, .$$ Using Lemma \ref{ancestral},
this is bounded up to some fixed
multiplicative constant by
$$\sum_{r=1}^{\lfloor n^{1/2+\eta}\rfloor}\sum_{h'=\lfloor n^\eta\rfloor}^{\lfloor n^{1/2+\eta}\rfloor}
\max_{n^\eta\leq h\leq h'}\wh{P}^{(x_r),h'}\left( \anc_T^V(c_T(v)\geq
  c\log n,h)\geq h^\gamma\right)\, .$$ By conditioning on the types of
the ancestors of $V$ under $\wh{P}^{(x_r),h'}$, we may again use (ii)
in Lemma \ref{sec:expon-contr-branch-1} to argue that the number of
ancestors $v$ of $V$ such that $c_F(v)\geq c\log n$ and $|v|\geq
|V|-h$ is a sum of independent Bernoulli random variables with
parameters of the form 
$$\frac{\zeta^{(i)}(|\bw|\geq c\log n)}{m_{ij}}\leq
n^{-c\eps}\frac{\max_{i\in
  [K]}\zeta^{(i)}(\exp(\eps|\bw|))}{\min_{ij:m_{ij\neq 0}}m_{ij}}\, ,$$ 
so
that it is stochastically bounded by a sum of $h$ i.i.d.\ Bernoulli
random variables with parameter $h^{-c\eps/(1/2+\eta)}\leq h^{-c\eps}$
as soon as $0<\eta<1/2$. Now if $X_1,\ldots,X_n$ are i.i.d.\ Bernoulli
random variables with common parameter $p$ on some probability space
$(\Omega,\a,P)$, $P(\sum_{i=1}^nX_i\geq x)\leq \exp(-\lambda
x+np(\exp(\lambda)-1))$ for every $\lambda>0$, by Markov's
inequality. Therefore, by choosing $\lambda$ small enough, for some
$C>0$, and by unconditioning on the sequence of types in the end,
$$\wh{P}^{(x_r),h'}\left( \anc_T^V(c_T(v)\geq
  c\log n,h)\geq 2h^{1-c\eps}\right)\leq \exp(-Ch^{1-c\eps})\, ,$$
hence the result by choosing $\gamma\in(1-c\eps,1)$, provided
$n^\eta\leq h$ and $n$ is large enough.

Point (iii) is obtained following similar lines. We are even going to
prove the result with $\anc^{u(k)}_F(c_F(v)=m,h)$ replaced by
$\anc^{u(k)}_F(c_F(v)=m,h)$.  Fix $\Gamma>0$, then the probability
under consideration is bounded up to a multiplicative constant by
$$\sum_{r=1}^{\lfloor n^{1/2+\eta}\rfloor}\sum_{h'=\lfloor n^\eta\rfloor}^{\lfloor n^{1/2+\eta}\rfloor}
\wh{P}^{(x_r),h'}\left(\max_{1\leq m\leq C\log n}\max_{n^\gamma\leq
    h\leq h'} m^Mh^{-1} \anc_T^V(c_T(v)\geq m,h)\geq
  C_M\right)+o(n^{-\Gamma})\, ,$$ for some $C>0$, where we used point
(i) of the lemma to bound the maximal degree of a vertex in the tree,
and then Lemma \ref{ancestral}. The first term is in turn bounded by
$$n^{1+2\eta}\max_{i\in[K]}\max_{n^\eta\leq h\leq h'\leq n^{1/2+\eta}}
\max_{1\leq m\leq C\log n}\wh{P}^{(i),h'}(m^Mh^{-1}
\anc_T^V(c_T(v)\geq m,h)\geq C_Mm^{-M}h)\, .$$ Using again the large
deviation inequality as above, we obtain that the probability on the
right hand-side is bounded by $\exp(-\lambda h(C_M m^{-M}-C'\exp(-\xi
m)))$ for some $\lambda>0$ and where $C',\xi$ are such that 
$$\max_{i\in[K]}
\zeta^{(i)}(|\bw|\geq x)\leq 2^{-1}C'\exp(-\xi x)\, .$$ By choosing
$C_M$ large enough, we obtain that
$$\exp(-\lambda h(C_M m^{-M}-C'\exp(-\xi
m)))\leq \exp(-\lambda' hm^{-M})$$ for some $\lambda'>0$ and every
$m\geq 1$. Since we consider only the terms $m\leq C\log n$ and $h\geq
n^\eta$, we get that $hm^{-M}\geq n^\eta(C\log n)^{-M}$ which is $\geq
n^{\eta/2}$ for large $n$, and this allows to conclude.  \cq

\subsection{Finite-dimensional marginals}\label{sec:finite-dimens-marg}

We are now able to show the convergence of finite-dimensional
marginals for the spatial process $S$. Notice the slightly loosened
hypotheses when compared with Theorem \ref{T2}. 

\begin{prp}\label{sec:finite-dimens-marg-1}
Assume that $\zzeta$ satisfies {\rm (H)} and admits small exponential
moments. Also assume that the spatial displacement laws $\nu_{i,\bw}$
are non-degenerate, centered and have a finite variance
$\langle\nu_{i,\bw},|y|_2^2\rangle$, satisfying 
\begin{equation}
  \max_{i\in[K]}\langle\nu_{i,\bw},|y|_2^2\rangle=O(|\bw|^D)\, ,
\end{equation}
for some $D>0$. Then, for any $\bx\in[K]^\N$, jointly with the
convergence {\rm (i)} in Theorem \ref{T1}, for any $0\leq
s_1<s_2<\ldots<s_k$, it holds that under $P^\bx$, 
$$\left(\frac{S^F_{\lfloor ns_j\rfloor}}{n^{1/4}},1\leq j\leq k\right)
\build\longrightarrow_{n\to\infty}^d\left(R_{s_j},1\leq j\leq
  k\right)\, ,$$ where $R$ has the law described in Theorem \ref{T2}.
If $\bx={\bf i}$, it also holds jointly with {\rm (ii)} in Theorem
\ref{T1}.
\end{prp}

\proof Since the proof is very similar to that of \cite[Proposition
28]{jfmgm05} and the general case does not involve more sophisticated
tools besides notational annoyments, we give the detailed argument
only for $k=2$.

Making use of Skorokhod's representation theorem, we will assume that
the discrete snakes we are considering are defined on a probability
space $(\Omega,\a,P)$ so that following holds. First, this space
supports a sequence of processes $H^n,\Lambda_i^n,i\in[K]$ with same
distribution as $n^{-1/2}H^F_{\lfloor n\cdot\rfloor},\Lambda^F_i(\lfloor n\cdot\rfloor),i\in[K]$
under $P^\bx$, such that $H^n$ converges almost-surely to a process
$B$ which is $2\sigma^{-1}$ times a standard Brownian motion, while
$\Lambda^n_i/a_i$ converges almost-surely to the identity
function. Otherwise said, the convergence of Theorem \ref{T1} is
almost-sure. Second, for every $n\geq 1$, the processes
$(n^{1/2}H^n_{k/n},\Lambda^n_i(k/n),i\in[K],k\geq 0)$ are the height
and type-counting processes of a unique multitype random forest
$(F^n,e^n)\in\FF^{(K)}$, and we assume that the discrete snakes are
defined by using a family of random variables $(Y^n_u,u\in F^n)$ that
are supported on $(\Omega,\a,P)$. That is, given $F^n$, the vectors
$(Y^n_{ul},1\leq l\leq c_{F^n}(u)), u\in F^n$ are independent with
respective distributions $\nu_{e^n(u),\bw_{F^n}(u)}$, and we let
$$S^n_k=\sum_{v\vdash u_{F^n}(k)}Y^n_v\, ,\qquad k\geq 0$$
be the associated discrete snake.

Let $s<t$ be given. For simplicity write
$u=u_{F^n}(\lfloor ns\rfloor),u'=u_{F^n}(\lfloor nt\rfloor)$ and $\check{u}$ the most recent
common ancestor of $u,u'$. The vector $(S^n_{ns},S^n_{nt})$ is the
image by the application $(x,y,z)\mapsto (x+z,y+z)$ of the vector
$({\sf L}^n_1,{\sf R}^n_1,{\sf C}^n)$ where
$${\sf L}^n_1=\sum_{\check{u}\vdash v\vdash u,v\neq \check{u}}Y^n_v \, ,\quad 
{\sf R}^n_1=\sum_{\check{u}\vdash v\vdash u',v\neq \check{u}}Y^n_v \,
,\quad {\sf C}^n=\sum_{v\vdash \check{u}}Y^n_v \, .$$ We aim at
proving that given $B$, $n^{-1/4}({\sf L}^n_1,{\sf R}^n_1,{\sf C}^n)$
converges in distribution to a triple of independent Gaussian
variables with respective variances $\Sigma^2 B_s,\Sigma^2
B_t,\Sigma^2 \check{B}_{s,t}$, where $\check{B}_{s,t}=\inf_{s\leq
  u\leq t}B_u$.

First note that conditionally on $F^n$, ${\sf L}^n_1,{\sf R}^n_1,{\sf
  C}^n$ are almost independent as they are sums involving terms which
are all independent but for one: if $l,l'$ are such that $u\in
\check{u}lF^n_{\check{ul}}$ and $u'\in\check{u}l'F^n_{\check{u}l'}$,
i.e.\ $u$ and $u'$ are in the $l$-th and $l'$-th subtree pending from
$\check{u}$, then $Y^n_{\check{u}l}$ and $Y^n_{\check{u}l'}$ may be
dependent. However, letting ${\sf L}^n_1=Y^n_{\check{u}l}+{\sf L}^n_2$
and ${\sf R}^n_1=Y^n_{\check{u}l'}+{\sf R}^n_2$, the vector $({\sf
  L}^n_2,{\sf R}^n_2,{\sf C}^n)$ has three independent components.  It
is thus sufficient to prove that $n^{-1/4}\max_l Y^n_{\check{u}l}$ and
$n^{-1/4}Y^n_{\check{u}}$ both converge to $0$ as $n\to\infty$, and to
show the individual convergences of $n^{-1/4}{\sf L}^n,n^{-1/4}{\sf
  R}^n,n^{-1/4}{\sf C}^n$ to three Gaussian variables with the correct
variances, where ${\sf L}^n$ and ${\sf R}^n$ are defined as ${\sf
  L}^n_1$ and ${\sf R}^n_1$, but allowing the term $v=\check{u}$ in
the sum.

Using Lemma \ref{sec:expon-contr-branch-2} for $\Gamma>1$ and applying
the Borel-Cantelli Lemma, we have that almost-surely, for $n$ large
enough, $c_{F^n}(u(k))\leq C\log n$ for every $k\leq nt$ and for some
$C>0$. Conditionally on $F^n$ and on this event, the expectation of
$\max_l (Y^n_{\check{u}l})^2$ is bounded above by $(C\log n)^D$ by our
assumption on moments, from which it follows that $n^{-1/4}\max_l
Y^n_{\check{u}l}$ does converge to $0$. The reasoning is similar for
estimating $n^{-1/4}Y^n_{\check{u}}$.

Let us now deal with ${\sf L}^n$, the case of ${\sf R}^n$ being similar.  The
sum defining ${\sf L}^n$ is a sum along a part of the ancestral line of $u$,
whose height is approximately $\sqrt{n}(B_s-\check{B}_{s,t})$. Note
that a.s., $B_s-\check{B}_{s,t}>0$ by standard properties of Brownian
motion, so that the sum in question has an order of $\sqrt{n}$
terms. Write
$${\sf L}^n=\sum_{\check{u}\vdash v\vdash u}Y^n_v
\ind_{\{c_{F^n}(u)\leq c\}}+\sum_{\check{u}\vdash v\vdash
  u}Y^n_v \ind_{\{c_{F^n}(u)> c\}}\, .$$
We call these two sums ${\sf L}^n_c$ and $\tilde{{\sf L}}^n_c$. 

We first argue that 
\begin{equation}\label{eq:18}
\lim_{c\to\infty}\limsup_{n\to\infty}
P(n^{-1/4}|\tilde{{\sf L}}^n_c|>\eps|B)=0\, .
\end{equation} 
Letting $h=|u|-|\check{u}|$, we
noticed that given $B$, $h$ is of the order $n^{1/2}$. In particular,
using the Borel-Cantelli Lemma and (iii) in Lemma
\ref{sec:expon-contr-branch-2} for $M=D+2$, for some constant $C_D$,
and every $m\geq 1$, a.s.\ $\anc_{F^n}^u(c_{F^n}(v)=m,h)\leq
C_Dm^{-D-2}h$, for every large enough $n$. Then by Chebychev's
Inequality and conditional independence of the terms,
\begin{eqnarray*}
  P(n^{-1/4}|\tilde{{\sf L}}^n_c|>\eps|B,F^n)&\leq & n^{-1/2}\eps^{-2}
  \sum_{i\in[K]} \sum_{\bw\in\W_K,|\bw|\geq c}\sum_{1\leq l\leq |\bw|}
  \anc_{F^n}^u(i,\bw,l,h)\langle\nu_{i,\bw},y_l^2\rangle\\
  &\leq & Cn^{-1/2}\eps^{-2} \sum_{m\geq c}m^D
  \anc_{F^n}^u(c_{F^n}(v)= m,h)\\
  &\leq & C_D\eps^{-2}n^{-1/2}h\sum_{m\geq c}m^{-2}\, ,
\end{eqnarray*}
at least for ($\omega$-dependent) large enough $n$.  Notice that this
upper bound does not depend on $F^n$. Since $n^{-1/2}h$ converges to
$B_s-\check{B}_{s,t}$, taking the conditional expectation given $B$
and applying the reverse Fatou Lemma, we obtain 
\begin{eqnarray*}
  \limsup_{n\to\infty}P(n^{-1/4}|\tilde{{\sf L}}^n_c|>\eps|B)
  &\leq &E[\limsup_{n\to\infty}P((n^{-1/4}|\tilde{{\sf L}}^n_c|>\eps|B,F^n)|B]\\
  &\leq & C_D\eps^{-2}(B_s-\check{B}_{s,t})\sum_{m\geq c}m^{-2}\, ,
\end{eqnarray*}
which goes to $0$ as $c\to\infty$. 

On the other hand, fixing $\eta,\eps>0$, we may apply Lemma
\ref{sec:proof-theorem-reft2-1} and Borel-Cantelli to obtain that
a.s., for large enough $n$, and every $0\leq k\leq nt,n^{\eta}\leq
h\leq |u(k)|, \bw\in\W_K$ with $|\bw|\leq c$ and $1\leq l\leq |\bw|$,
\begin{equation}\label{eq:19}
|\anc_{F^n}^{u(k)}(i,\bw,l,h)-ha_ib_{w_l}\zeta^{(i)}(\bw)|\leq
h^{1/2+\eta}\leq \eps ha_ib_{w_l}\zeta^{(i)}(\bw)\, .
\end{equation} 
Now, we can redisplay ${\sf L}^n_c$ as the finite sum
$$\sum_{i\in[K]}\sum_{\bw\in \W_K,|\bw|\leq c}\sum_{l=1}^{|\bw|}
\sum_{\check{u}\vdash v\vdash u}Y^n_v\ind_{\{e(v)=i,\bw_{F^n}(v)=\bw,
  u\in vlF^n_{vl}\}}\, ,$$ and given $F^n$, the last summation has
$\anc_{F^n}^u(i,\bw,l,h)$ i.i.d.\ terms (for $h=|u|-|\check{u}|$) with
variance $\langle\nu_{i,\bw},y_l^2\rangle$. Using (\ref{eq:19}), we
can apply the the central limit theorem conditionally on $B,F^n$, and
obtain that $n^{-1/4}{\sf L}^n_c$ converges in distribution to a
Gaussian variable with variance $(B_s-\check{B}_{s,t})\Sigma_c^2$
where
$$\Sigma_c^2=\sum_{i\in[K]}a_i\sum_{\bw\in\W_K,|\bw|\leq
  c}\zeta^{(i)}(\bw)\sum_{l=1}^{|\bw|}b_{w_l}\langle\nu_{i,\bw},y_l^2\rangle\,
.$$ As $c\to\infty$, this converges to the constant $\Sigma^2$ of the
statement of the theorem, which implies the result when combined with
(\ref{eq:18}).

To complete the proof, it remains to show the convergence of
$n^{-1/4}{\sf C}^n$. The argument is essentially the same, where the height
$h$ under consideration is of the order $\sqrt{n}\check{B}_{s,t}$
rather than $\sqrt{n}(B_s-\check{B}_{s,t})$. The argument above is
unchanged in the case $\check{B}_{s,t}>0$, but is not valid anymore if
$\check{B}_{s,t}=0$. In that case, we claim that the result is in fact
trivial as we have $|\check{u}|=0$, so that ${\sf C}^n=0$. Indeed, assume
for a moment that $\bx={\bf i}$ for some $i\in[K]$.  By Theorem \ref{T1},
and plugging a new random element in our use of Skorokhod's
representation Theorem, we may assume that
$(n^{-1/2}\Upsilon^{F^n}_{\lfloor ns\rfloor},s\geq 0)$ converges a.s.\ to a
multiple of the local time $L^0$ of $B$ at level $0$. Since
$\check{B}_{s,t}=0$ and $B_sB_t>0$, there is an increase time of $L^0$
between times $s$ and $t$. Therefore, the function $\Upsilon^{F^n}$
increases as well between these times, at least for large enough
$n$. This means that $u$ and $u'$ are in two different tree components
of $F^n$, and thus $\check{u}=\varnothing$.

The case where $\bx$ is any element of $[K]^\N$ is a slight
elaboration of the preceding argument, which we briefly sketch. If
$\check{B}_{s,t}=0$ but $|\check{u}|>0$ for infinitely many $n$, this
means that for these values of $n$, the vertices $u,u'$ belong to the
same tree component of $F^n$, whose root is of type $i$, say. Now,
$n^{-1/2}H^{\Pi^{(i)}(F^n)}_{\lfloor n\cdot\rfloor}$ converges to a multiple of the
Brownian motion $B_{\cdot/a_i}$ by the proof of Theorem \ref{T1},
while $n^{-1/2}\Upsilon^{\Pi^{(i)}(F^n)}_{\lfloor n\cdot\rfloor}$ converges to a
multiple of the local time of the latter. Necessarily, there is a time
of increase of this local time between the times $a_is$ and $a_it$,
hence the function $n^{-1/2}\Upsilon^{\Pi^{(i)}(F^n)}_{\lfloor n\cdot\rfloor}$
increases also during that time interval, whose length corresponds
asymptotically to the fraction of vertices of type $i$ that appear
between $u$ and $u'$ in depth-first order. This is a contradiction
with the fact that $u$ and $u'$ belong to the same tree component of
$F^n$. \cq

\subsection{H\"older norm bounds}\label{sec:holder-norm-bounds}

\begin{prp}\label{sec:proof-theorem-reft2}
  Let $\mmu$ satisfy the basic assumption {\rm (H)}, and admit small
  exponential moments.  For every $A>0$ and $\alpha\in(0,1/4)$, for
  every $\eps>0$, there exists $C>0$ such that for every
  $\bx\in[K]^\N$,
$$\sup_{n\in\N}P^{\bx}\left(\sup_{s\neq t\in[0,A],ns,nt\in\Z_+}
  \frac{|H^F_{ns}-H^F_{nt}|}{\sqrt{n}|s-t|^\alpha}>C\right)\leq
\eps\, .$$ 
\end{prp}

\rem As will be shown below, in the particular case $K=1$, we are able
prove the same assertion with $\alpha<1/2$ rather than
$\alpha<1/4$. In \cite{jfmgm05}, we could also obtain the result for
$\alpha<1/2$ because of the particular nature of the multitype trees
we were considering, i.e.\ alternating types, corresponding to an
antidiagonal mean matrix. In the general case, our method of
approximation by the monotype case does not seem to be fine enough to
obtain the best estimate. 

\bigskip

\proof We first prove the result in the case $K=1$, and let
$\mu:=\mu^{(1)}$, $P:=P^{{\bf 1}}$. Our proof is partly inspired from
that of \cite[Theorem 1.4.4]{duqleg02}, of which we can interpret
Proposition \ref{sec:proof-theorem-reft2} to be a discrete
counterpart.

Recall e.g.\ from \cite{duqleg02} that if we let
$$V^{\bf f}_n= \sum_{k=1}^n(c_{\bf f}(u(k))-1)\, , \qquad n\geq 0$$
be the {\em \L ukaciewicz} walk associated with the forest $\ff$,
then the height process of $\ff$ is given by
\begin{equation}\label{eq:1}
  H^\ff_n=\#\left\{k\in\{0,1,\ldots,n-1\}:V^\ff_k=\min_{k\leq l\leq n}
    V^\ff_l\right\}.
\end{equation}
Under $P$, $V^F$ is a random walk on $\Z$ with centered step
distribution $\mu(\cdot+1)$ on $\{-1,0,1,2,\ldots\}$. 

Now, suppose $0\leq s<t\leq A$ are such that $ns,nt\in\Z_+$. Write
$\lambda(x)=\max\{ l\in[0, ns] :V^F_l\leq x\}$.  Using
(\ref{eq:1}), we have
\begin{eqnarray}
  H^F_{nt}-H^F_{ns}&=&
  \#\left\{k\in[ns,nt):V^F_k=\min_{k\leq l\leq nt}V^F_l
  \right\}\nonumber\\
  & &-\#\left\{\lambda\left(\min_{ns\leq l\leq nt}V^F_l\right)<
    k<ns: V^F_l=\min_{k\leq l\leq ns}V^F_l\right\}\, , \label{eq:2}
\end{eqnarray}
and the rest of the proof will consist in estimating the moments of
the two terms above, which correspond to the lengths of the branches
of $F$ from $u(ns),u(nt)$ down to their most recent common ancestor.
By the time reversal property for walks,
$$(\wh{V}^{(n)}_k=V_n^F-V_{n-k}^F,0\leq k\leq n)
\build=_{}^d(V^F_k,0\leq k\leq n),$$ the first term in (\ref{eq:2}) is
equal in distribution, under $P$, to
$$W_{n(t-s)}:=\#\left\{1\leq k\leq n(t-s):V^{F}_k=
\max_{0\leq l\leq k} V^{F}_l\right\},$$ 
the number of (weak) records of $V^F$ before time $n(t-s)$. 

Let $M_{n}=\max_{0\leq k\leq n}V^F_k$. Let $\tau_0=0$, and
$\tau_i,i\geq 1$ be the $i$-th record time, i.e.\ the $i$-th time
$\tau\geq 1$ such that $V^F_{\tau}=M_{\tau}$. Then it is easy and
well-known that $(\tau_i-\tau_{i-1},i\geq 1)$ are a sequence of
i.i.d.\ random variables.  Moreover, since $V^F$ is centered and
its increments have finite second moment under $P$, it
is a consequence of the proof of \cite[XII,7 Theorem 1a]{feller2} and
the discussion before that the Laplace exponent $\phi(s)=-\log
E[\exp(-s\tau_1)]\sim Cs^{1/2}$ as $s\to0$ for some
$C>0$ (Feller considers the case of strict ladder epochs, but the
treatment of weak ones is similar). Now, for any $p>1$, and integer
$u$,
\begin{eqnarray*}
  E[W^p_u]&=& p\int_0^{\infty}x^{p-1} P
  \left(W_u\geq x\right)\d x\\
  &=& p\int_0^{\infty}x^{p-1} P
  \left(\sum_{i=1}^x(\tau_i-\tau_{i-1})\leq u\right)\d x\\
  &\leq & pe \int_0^{\infty}x^{p-1} E
  \left[\exp\left(-\sum_{i=1}^x\frac{\tau_i-\tau_{i-1}}{u}\right)\right]\d x
  \leq C' \phi(u^{-1})^{-p}\leq C'' u^{p/2},
\end{eqnarray*}
for some $C',C''>0$ and every $u$.  
Since $ns,nt$ are distinct integers, we showed that
$E[W_{n(t-s)}^p]\leq C_1 n^{p/2}|s-t|^{p/2}$ uniformly in such
$n,s,t$, where $C_1=C_1(\mu,p)>0$.

Let us now handle the second term in (\ref{eq:2}). Using
time-reversal, we see that this equals
$$\#\left\{n(t-s)< k < nt\wedge
  \kappa\left(\max_{1\leq l\leq n(t-s)}V^{F}_l\right):
  V^{F}_k=\max_{n(t-s)\leq l\leq k}V^{F}_l \right\}$$ in distribution,
where $\kappa(x)=\min\{k\geq n(t-s):V^{F}_k\geq x\}$ (with the
convention $\min\varnothing=\infty$).  By using Markov's property at
time $n(t-s)$, this has same distribution as $W_{ns\wedge
  \kappa(\wt{M}_{n(t-s)}-\wt{V}_{n(t-s)})-1}$, where $W$ is defined as
above, while $\wt{V}$ is an independent copy of $V^{F}$ with maximum
process $\wt{M}$.  By monotonicity this is less than
$W_{\kappa(\wt{M}_{n(t-s)}-\wt{V}_{n(t-s)})}$.  Let us prove that
$E[W_{\kappa(x)}^p]\leq Cx^p$ for every $x\geq 0$, for some $C>0$.  To
this end, notice that
\begin{equation}\label{eq:3}
M_{n(t-s)}=\sum_{i=1}^{W_{n(t-s)}}(V^{F}_{\tau_i}-V^{F}_{\tau_{i-1}}),
\end{equation}
and it is a classical result of fluctuation theory that the variables
$V^{F}_{\tau_i}-V^{F}_{\tau_{i-1}}$ are independent with common
distribution
$P(V^{F}_{\tau_1}=i)=\mu([i+1,\infty)),i\geq 0$, so
their mean is $\sigma^2/2$, where $\sigma^2$ is the variance of $\mu$,
and notice that these variables have small exponential moments. Now,
the usual large deviations theorem shows that for some $a,N>0$ and for
every $n\geq N$,
\begin{equation}\label{eq:4}
  P\left(\frac{\sum_{i=1}^n(V^{F}_{\tau_i}-
      V^{F}_{\tau_{i-1}})}{n}< \frac{\sigma^2}{4}
  \right)\leq \exp(-an).
\end{equation}
Now, using (\ref{eq:3}) in the second equality,
\begin{eqnarray*}
  E[W_{\kappa(x)}^p] &=& p\int_0^{\infty}u^{p-1}
  P(W_{\kappa(x)}> u)\d u\\
  &=& p\int_0^{\infty}u^{p-1} P\left(\sum_{i=1}^{
      \lceil u \rceil}(V^{F}_{\tau_i}-V^{F}_{\tau_{i-1}})< x\right)\d u\\
  &=& px^p\int_0^{\infty}v^{p-1} \, P\left(\sum_{i=1}^{
      \lceil xv \rceil}(V^{F}_{\tau_i}-V^{F}_{\tau_{i-1}})< 
    \frac{xv}{v}\right)\d v\\
  &\leq & px^p\left(\int_0^{4\sigma^{-2}}v^{p-1}\d v+
    \int_{4\sigma^{-2}}^{\infty}
    v^{p-1}P\left(\sum_{i=1}^{
        \lceil xv \rceil}(V^{F}_{\tau_i}-V^{F}_{\tau_{i-1}})< 
      \frac{\sigma^2 xv}{4}\right)\d v\right).
\end{eqnarray*}
Now, as soon as $x$ is large enough, i.e.\ $4x\sigma^{-2}\geq N$,
where $N$ is defined before (\ref{eq:4}), the probability in the
second integral is bounded by $\exp(-axv)\leq \exp(-v)$ if we further
ask $x>a^{-1}$. Thus the wanted bound on $E[W_{\kappa(x)}^p]$. By the
independence of $\wt{V}$ and $V^F$, we conclude that
\begin{eqnarray}
\label{gnst}
E[(W_{\kappa(\wt{M}_{n(t-s)}-\wt{V}_{n(t-s)})})^p]
&\leq& CE[(M_{n(t-s)}-V^{F}_{n(t-s)})^p] \\
&\leq & 2^{p-1}C\left(1+\left(\frac{p}{p-1}\right)^p\right)
E[(V^{F}_{n(t-s)})^p] , \nonumber
\end{eqnarray}
where we used Doob's inequality $E[M_{n(t-s)}^p]\leq (p/(p-1))^p
E[(V^{F}_{n(t-s)})^p]$, since $V^{F}$ is centered. Now we use the
following consequence of Rosenthal's inequality \cite[Theorem
2.10]{petrov95}: if $X_1,\ldots,X_n$ are independent centered random
variables (not necessarily identically distributed) defined on some
probability space $(\tilde{\Omega},\tilde{\a},\tilde{P})$, then for
every $p\geq 2$ there exists $C(p)$ such that
\begin{equation}\label{rosenthal}
  \tilde{E}[|X_1+\ldots+X_n|^p]\leq C(p)n^{p/2-1}\sum_{i=1}^n
  \tilde{E}[|X_i|^p]. 
\end{equation}
This shows that $E[(V^{F}_{n(t-s)})^p]\leq
C'(p)n^{p/2} |s-t|^{p/2}$ for some $C'(p)>0$, for every $s,t$ such
that $ns,nt\in\Z_+$, and therefore the same kind of upper bound holds
for the quantity in (\ref{gnst}).

Putting things together, we have obtained that for every $p\geq 2$ and
some $C_2=C_2(\mu,p)>0$,
$$\sup_{n\geq 1}\sup_{s,t\geq 0, ns,nt\in \Z_+}
E\left[\left|\frac{H^{F}_{ns}-H^{F}_{nt}}{ \sqrt{n}}\right|^p \right]
\leq C_2|s-t|^{p/2}.$$ Let now $H^{F}_{\{s\}},s\geq 0$ be defined by
linear interpolation between linear abscissa. Then, it is elementary
that 
$$\sup_{n\geq 1}\sup_{s\neq t\geq 0}
E\left[\left|\frac{H^{F}_{\{ns\}}-H^{F}_{\{nt\}}}{ \sqrt{n}}\right|^p
\right] \leq C_3|s-t|^{p/2}\, ,$$ for some $C_3>0$. The uniform
estimate in Kolmogorov's criterion \cite[Theorem 3.4.16]{Stroock93}
finally entails the result, for $K=1$.

The general case is obtained by using the contraction function
$\Pi^{(1)}$. We have, for $0\leq s,t\leq A$ and $ns,nt\in\Z_+$,
$$
|H^F_{ns}-H^F_{nt}|\leq \frac{|H^{\Pi^{(1)}(F)}_{\Lambda^F_1(ns)-1}-
  H^{\Pi^{(1)}(F)}_{\Lambda^F_1(nt)-1}|}{a_1b_1} +2\max_{0\leq k\leq
  An}\left|H^F_k-
  \frac{H^{\Pi^{(1)}(F)}_{\Lambda_1^F(k)-1}}{a_1b_1}\right|\, .
$$
On the one hand, by the case $K=1$ and Proposition \ref{monotype}, we
have with high $P^\bx$-probability, uniformly in $n\in\N$,
$$\frac{|H^{\Pi^{(1)}(F)}_{\Lambda^F_1(ns)-1}-
  H^{\Pi^{(1)}(F)}_{\Lambda^F_1(nt)-1}|}{\sqrt{n}|s-t|^\alpha}\leq
C\frac{|n^{-1}\Lambda^F_1(ns)-n^{-1}\Lambda^F_1(nt)|^{\alpha}}{
|s-t|^{\alpha}}\leq C\, ,$$
since $\Lambda^F$ is a counting function. On the other hand, using the
same inequalities as around (\ref{eq:8}), 
\begin{eqnarray*}\lefteqn{\max_{0\leq k\leq An}
\left|H^F_k-
  \frac{H^{\Pi^{(1)}(F)}_{\Lambda_1^F(k)-1}}{a_1b_1}\right|}\\
&\leq&
\max_{0\leq k\leq An}\left|H^F_k-
  \frac{\anc_F(1,u(k))}{a_1b_1}\right|+\max_{0\leq k\leq An}
\frac{|H^{\Pi^{(1)}(F)}_{k-1}-H^{\Pi^{(1)}(F)}_k|+1}{a_1b_1}\, .
\end{eqnarray*}
According to Proposition \ref{HEIGHT}, the first term on the
right-hand side is bounded above by $n^{1/4+\gamma}$ with high
probability for large $n$, where we chose $0<\gamma<1/4-\alpha$, so
that $n^{1/4+\gamma}\leq \sqrt{n}|s-t|^{\alpha}$ for every large $n$
(recall that $ns,nt\in\Z_+$ so that $|s-t|\geq n^{-1}$). As for the
second term, by \cite[Lemma 21]{jfmgm05} (while the statement is on
conditioned trees, the first part of its proof yields the result on
forests), it holds that it is $o(n^{\gamma})$ in $P^\bx$-probability
for any $\gamma>0$, which gives the wanted bound. \cq

\subsection{Tightness}\label{sec:tightness}

This section is devoted to the proof of last building block needed to
prove Theorem \ref{T2}, namely

\begin{prp}\label{sec:tightness-1}
  Under the hypotheses of Theorem \ref{T2}, for every $\bx\in[K]^\N$,
  the laws of the processes $((n^{-1/4}S^F_{\lfloor ns\rfloor},s\geq 0),n\geq 1)$
  under $\P^\bx$ are tight in the Skorokhod space $\mathbb{D}(\R+,\R)$.
\end{prp}

\proof In this proof, $C_1,\ldots,C_{10}$ will be denoting strictly
positive constants. Our first task is to obtain an upper bound for
expectations of the form $\E^\bx[|S^F_{nt}-S^F_{ns}|^p]$.
To this end, we first choose $\xi,D$ so that
the moment condition (\ref{eq:5}) holds, and write $p=8+\xi$. Also,
fix $0<\alpha<1/4$ so that $\alpha(8+\xi)>2$, $0<\eta<1/4$, and $c>0$
such that $c<\eta/(4\log K)$. According to Proposition
\ref{sec:holder-norm-bounds}, we may choose $C_1>0$ so that if
$$A_n=\left\{\max_{0\leq s\neq t\leq A:ns,nt\in\Z_+}
  \frac{|H^F_{ns}-H^F_{nt}|}{n^{1/2}|s-t|^\alpha}\leq C_1\right\}\,
,$$ then $P^{\bx}(A_n)\geq 1-\eps/2$ for every $n\geq 1$. We let $B_n$
be the intersection of $A_n$ with the three events $\{\max_{0\leq
  k\leq An}c_F(u(k))\leq C\log n\}$,
$$\left\{\max_{0\leq k\leq An}\max_{n^\eta\leq h\leq
    |u(k)|}\max_{\bw\in \W_K,1\leq l\leq |\bw|}\max_{i\in [K]}\frac{
    |\anc^{u(k)}_F(i,\bw,l,h)-ha_i
    b_{w_l}\zeta^{(i)}(\bw)|}{h^{1/2+\eta}}\leq 1\right\}\, $$ and
$$\left\{\max_{0\leq k\leq An}\max_{n^\eta\leq h\leq |u(k)|}
  \anc_F^{u(k)}(c_F(v)\geq c\log n,h)\leq h^\gamma\right\}\, ,$$ where
$C>0,0<\gamma<1$ are chosen so that $P^\bx(B_n)\geq 1-\eps$ for every
$n$ sufficiently large, which is possible according to Lemmas
\ref{sec:proof-theorem-reft2-1} and \ref{sec:expon-contr-branch-2}.
We take $n\in\N$, choose $0\leq s\neq t\leq A$ be such that
$k=ns\in\Z_+,k'=nt\in\Z_+$, and write $u=u(k),u'=u(k')$.  Then, by
definition, we have
$$S^F_k-S^F_{k'}=\sum_{v\vdash u,|v|>|\check{u}|}Y_v-
\sum_{v\vdash u',|v|>|\check{u}|}Y_v\, ,$$ whenever $\check{u}$ is the
most recent common ancestor to $u$ and $u'$.  Assume
$u=\check{u}rw,u'=\check{u}r'w'$ for some $r\neq r'\in\N$ and
$w,w'\in{\cal U}$. This allows to redisplay the previous expression as
\begin{eqnarray*}
  S^F_k-S^F_{k'}&=& (Y_{\check{u}r}-Y_{\check{u}r'})\\
  & & +\sum_{\bw\in\W_K}\sum_{1\leq l\leq |\bw|}\sum_{i\in[K]}
  \sum_{v\vdash u,v\neq u,v\in
    F^{(i)}}Y_{vl}\ind_{\{|v|>|\check{u}|,\bw_F(v)=\bw,
    vl\vdash u\}}\\
  & &-\sum_{\bw\in\W_K}\sum_{1\leq l\leq |\bw|}\sum_{i\in[K]}
  \sum_{v\vdash u',v\neq u',v\in
    F^{(i)}}Y_{vl}\ind_{\{|v|>|\check{u}|,\bw_F(v)=\bw,vl\vdash
    u'\}}\, .
\end{eqnarray*}
By construction, under $\E^\bx$, all the terms of this sum are
independent of each other conditionally on $F$, except possibly for
$Y_{\check{u}r}$ and $Y_{\check{u}r'}$. Let $h=|u|-|\check{u}|-1$ and
$h'=|u'|-|\check{u}|-1$, and let
$R(k,k')=|u|+|u'|-2|\check{u}|=h+h'+2$ be the number of random
variables of the form $Y_v$ that are involved in the expression
$S^F_k-S^F_{k'}$. Using (\ref{rosenthal}) for $p=8+\xi$ gives,
\begin{eqnarray}
  \lefteqn{\E^\bx\left[\left|S^F_k-S^F_{k'}\right|^p\, |\, F\right]}\\
  &\leq & C_2 R(k,k')^{p/2-1}\times
  \left(\begin{array}{l}
      \E^\bx[|Y_{\check{u}r}-Y_{\check{u}r'}|^p\, | \, F]\\
      +\sum_{\bw\in\W_K}\sum_{l=1}^{|\bw|}\sum_{i\in[K]}
      \anc_F^u(i,\bw,l,h)\langle\nu_{i,\bw},|y_l|^p\rangle\\
      +\sum_{\bw\in\W_K}\sum_{l=1}^{|\bw|}\sum_{i\in[K]}
      \anc_F^{u'}(i,\bw,l,h')\langle\nu_{i,\bw},|y_l|^p\rangle
\end{array}\right)\nonumber\\
&\leq &C_3 R(k,k')^{p/2-1}\times\left(\begin{array}{l}
    c_F(\check{u})^D\\
    +\sum_{m\geq 1}m^D\sum_{\bw\in \W_K,|\bw|=m}\sum_{l=1}^m
    \sum_{i\in[K]}
    \anc_F^u(i,\bw,l,h)\\
    +\sum_{m\geq 1}m^D\sum_{\bw\in \W_K,|\bw|=m}
    \sum_{l=1}^m \sum_{i\in[K]}
    \anc_F^{u'}(i,\bw,l,h')
\end{array}\right)\nonumber
\end{eqnarray}
Now, on the event $B_n$, we have $c_F(\check{u})^D\leq (C\log n)^D\leq
n^{1/2}|s-t|^\alpha$, since $|s-t|\geq 1/n$. 

It remains to bound the
two above sums, by symmetry it suffices to deal with the first one. In
the case where $h\leq n^\eta$ and on the event $B_n$, notice that
$$\sum_{m\geq 1}m^D\sum_{\bw\in \W_K,|\bw|=m}\sum_{l=1}^m
\sum_{i\in[K]} \anc_F^u(i,\bw,l,h)=\sum_{1\leq m\leq C\log
  n}m^D\anc_F^u(c_F(v)=m,h)\, ,$$ which is less than $(C\log n)^Dh\leq
(C\log n)^Dn^{\eta}$, and this in turn is less than
$n^{1/2}|s-t|^{\alpha}\geq n^{1/2-\alpha}$ since
$\eta<1/4<1/2-\alpha$. 

Assume now that $h\geq n^\eta$. Still on $B_n$, it holds that
$\anc_F^u(i,\bw,l,h)\leq ha_ib_{w_l}\zeta^{(i)}(\bw)+h^{1/2+\eta}$. We
split the sum under consideration into
$$
\sum_{1\leq m\leq c\log n}\!\!\!  m^D\!\!\!\sum_{\bw\in
  \W_K,|\bw|=m}\sum_{l=1}^m \sum_{i\in[K]}\anc_F^u(i,\bw,l,h)
+\!\!\!\sum_{c\log n< m\leq C\log n} m^D\anc_F^u(c_F(v)=m,h)
$$
The second term is bounded by $(C\log n)^Dh^{\gamma}$ by definition of
$B_n$ and since $h\geq n^\eta$, hence by $C_4h$. The
first term is bounded above by
$$C_5h\max_{i\in[K]}\sum_{\bw\in\W_K}|\bw|^{D+1}\zeta^{(i)}(\bw)
+C_5h^{1/2+\eta}(C \log n)^{D+1} \#\{\bw\in\W_K:|\bw|\leq c\log n\}\,
,$$ and since $\#\{\bw\in\W_K:|\bw|\leq c\log n\}\leq K^{c\log
  n}=n^{c\log K}\leq h^{c\log K/\xi}$ and by our choice of $c$, the
whole is bounded by $C_6h\leq C_6 R(k,k')$.

As in \cite[Proposition 27]{jfmgm05}, we argue that the number
$R(k,k')$, which equals $|u|+|u'|-2|\check{u}|=h+h'+2$, satisfies
$R(k,k')\leq C_7n^{1/2}|s-t|^{\alpha}$ on $B_n$, for all choices of
$n,s,t$. Hence by inspection of all the cases discussed above, we have
on $B_n$,
$$\E^\bx[|S^F_k-S^F_{k'}|^p|F]\leq C_8(n^{1/2}|t-s|^{\alpha})^{p/2}\,
,$$
so that 
$$\E^\bx\left[\left|\frac{S^F_{ns}-S^F_{nt}}{n^{1/4}}\right|^p\Big| B_n\right]
\leq C_9|t-s|^{\alpha p/2}\, .$$ As in the previous proof, if we let
$(S^F_{\{ns\}},s\geq 0)$ be the linearly interpolated version of
$(S^F_{\lfloor ns\rfloor},s\geq 0)$ between abscissa points of the form $k/n,k\geq
0$, then it is elementary that a similar bound holds up to taking a
larger $C_9$, this time for all $0\leq s\neq t\leq A$. Also, by our
choice of $\alpha$, we have $\alpha p/2>1$. Hence, an application of
Kolmogorov's criterion to $(n^{-1/4}S^F_{\{ns\}},s\geq 0)$ gives that
for every $A,\eps>0$, there exists $C_{10},\beta>0$ such that for
every $n\geq 1$,
\begin{equation}\label{eq:20}
  \P^\bx\left(\max_{0\leq t\neq s\leq A}
    \frac{|S^F_{\{nt\}}-S^F_{\{ns\}}|}{n^{1/4}|t-s|^{\beta}}\geq C_{10}
    \Big| B_n\right)\leq \eps\, , 
\end{equation} 
and since $P^{\bx}(B_n)\geq 1-\eps$, we may as well forget the
conditioning on $B_n$. 

This is enough to conclude that the laws of the continuous processes
$(n^{-1/4}S^F_{\{ns\}},s\geq 0)$ under $\P^\bx$ form a tight sequence,
and the result will follow from the fact that these processes are
respectively uniformly close to $(n^{-1/4}S^F_{\lfloor ns\rfloor},s\geq 0)$ over
compact intervals. This immediately comes from $\max_{1\leq k\leq
  An}|Y^F_{u(k)}|=o(n^{1/4})$ in probability, which itself can be
inferred from the fact that with high probability, no vertex $u(k)$
with $0\leq k\leq An$ has more than $C\log n$ children, and the moment
control (\ref{eq:5}). Hence the result. \cq

\subsection{Conditioned results: Theorem
  \ref{sec:conv-mult-snak}}\label{sec:cond-results:-theor}

Obtaining Theorem \ref{sec:conv-mult-snak} from Theorem
\ref{sec:convergence-gw-trees} now consists in reproducing faithfully
the proofs of Propositions \ref{sec:finite-dimens-marg-1} and
\ref{sec:tightness-1}, so we only sketch the plan of the proof. The
important results that are needed are generalizations to conditioned
measures of Lemmas \ref{sec:proof-theorem-reft2-1} and
\ref{sec:expon-contr-branch-2}, where $P^{\bx}$ must be replaced by
$P^{(i)}(\cdot|\#T^{(j)}=n)$. This is straightforward from the latter
lemmas and Lemma \ref{sec:cond-results:-theor-1}, which we use as we
did around (\ref{eq:22}). This is enough to obtain the exact analog of
Proposition \ref{sec:finite-dimens-marg-1} for the conditioned
probabilities $P^{(i)}(\cdot|\#T^{(j)}=n)$ (it is even simpler as the
issue encountered at the very end of the proof of that proposition
disappears).

The analog of the tightness statement (Proposition
\ref{sec:tightness-1}) is then a consequence of the following version
of Proposition \ref{sec:proof-theorem-reft2} for conditioned measures.

\begin{prp}\label{sec:cond-results:-theor-5}
Assume that $\mmu$ satisfies {\rm (H)} and admits small exponential
moments. 
  For every $i,j\in [K]$, for every $\alpha\in(0,1/4)$ and $\eps>0$,
  there exists $C>0$ such that
$$\sup_{n\in\N}P^{(i)}\left(\sup_{s\neq t\in[0,1],\#T s,\#T t\in\Z_+}
  \frac{|H^T_{\lfloor \#T s\rfloor}-H^T_{\lfloor \#T
      t\rfloor}|}{\sqrt{n}|s-t|^\alpha}>C\Big|\#T^{(j)}=n \right)\leq \eps\,
.$$
\end{prp}

\proof We rest on \cite[Theorem 24]{jfmgm05}, which is essentially the
monotype result ($K=1$), with the extra freedom that we consider the
conditioned law $P^{{\bf 1}}_r(\cdot|\#F=n)$ of a forest with $r$
components. Indeed, when applying the mapping $\Pi^{(j)}$ to $T$ under
$P^{(i)}(\cdot|\#T^{(j)}=n)$, one obtains such a conditioned forest
with a random number of roots, although the laws of these random
numbers form a tight sequence by Lemma
\ref{sec:cond-results:-theor-4}. Hence, up to conditioning, we can
assume that this number of roots is fixed and apply the monotype
result. The conclusion is then the exact analog of the last lines
of the proof of Proposition \ref{sec:proof-theorem-reft2}. Details are
left to the interested reader. \cq

\def\polhk#1{\setbox0=\hbox{#1}{\ooalign{\hidewidth
  \lower1.5ex\hbox{`}\hidewidth\crcr\unhbox0}}}

\end{document}